\documentclass{gtart}

\def\ifplaintex{\expandafter\ifx\csname documentclass\endcsname\relax}

\ifplaintex 
\else
\fi

\expandafter\ifx\csname beginpicture\endcsname\relax
\expandafter\ifx\csname documentclass\endcsname\relax
\input pictex \else
\input prepictex \input pictex \input postpictex \fi\fi

\def\gt{{\mathsurround=0pt\it $\cal G\mskip-2mu$eometry \&\ 
$\cal T\!\!$opology}}        

\def\gtp{{\mathsurround=0pt\it $\cal G\mskip-2mu$eometry \&\ 
$\cal T\!\!$opology $\cal P\!$ublications}}  

\def\lognumber#1{\def\thelognumber{#1}}
\def\volumenumber#1{\def\thevolumenumber{#1}}
\def\papernumber#1{\def\thepapernumber{#1}}
\def\volumeyear#1{\def\thevolumeyear{#1}}

\def\pagenumbers#1#2{\def\startpage{#1}\def\finishpage{#2}}
\def\published#1{\def\publishdate{#1}}
\def\proposed#1{\def\theproposer{#1}}
\def\seconded#1{\def\theseconders{#1}}
\def\received#1{\def\receiveddate{#1}}

\def\accepted#1{\def\accepteddate{#1}}

\long\def\asciiabstract#1{\long\def\theasciiabstract{#1}}

\let\\\par\let\thelognumber\relax
\let\thevolumenumber\relax\let\thepapernumber\relax
\let\thevolumeyear\relax\let\thesamplenumber\relax\let\startpage\relax
\let\finishpage\relax\let\publishdate\relax\let\receiveddate\relax
\let\reviseddate\relax\let\accepteddate\relax\let\theasciititle\relax
\let\theasciiauthors\relax
\let\theasciiabstract\relax
\let\theasciiemail\relax\let\theshortauthors\relax\let\theshorttitle\relax

\long\def\maketitlep{   

\count0=\startpage

\gt\hfill      
\beginpicture
\setcoordinatesystem units <0.33truein, 0.33truein> point at 2.2 0.9
\setplotsymbol ({$\cal G$})
\plotsymbolspacing=9truept
\circulararc 315 degrees from 0 1 center at 0 0
\setplotsymbol ({$\cal T$})
\circulararc 315 degrees from 1 -1 center at 1 0
\endpicture
\break
{\small\ifx\thesamplenumber\relax 
Volume \else Sample
\fi\thevolumenumber\ (\thevolumeyear)
\startpage--\finishpage\nl
Published: \publishdate}
\vglue 0.5truein plus 0.4fil minus 0.1truein

{\parskip=0pt\leftskip 0pt plus 1fil\def\\{\par\smallskip}{\ifplaintex\large
\else\Large\fi\bf\thetitle}\par\medskip}   

\vglue 0pt plus 0.1fil 

{\parskip=0pt\leftskip 0pt plus 1fil\def\\{\par}{\sc\theauthors}
\par\medskip}

\vglue 0pt plus 0.1fil 

{\small\parskip=0pt\let\newline\\
{\leftskip 0pt plus 1fil\def\\{\par}{\sl\theaddress}\par}
\expandafter\ifx\theemail\relax    
\relax\else\vglue 5pt plus 0.02fil minus 2pt\def\\{\stdspace{\rm 
and}\stdspace} 
\cl{Email:\stdspace\tt\theemail}\fi
\ifx\theurl\relax                  
\relax\else\vglue 5pt plus 0.02fil minus 2pt\def\\{\stdspace{\rm 
and}\stdspace}
\cl{URL:\stdspace\tt\theurl}\fi\par}

\vglue 7pt plus 0.3fil minus 3pt

{\bf Abstract}
\vglue 5pt plus 0.1fil minus 2pt

\theabstract

\vglue 7pt plus 0.3fil minus 3pt

{\bf AMS Classification numbers}\quad Primary:\quad \theprimaryclass

Secondary:\quad \thesecondaryclass

\vglue 5pt plus 0.3fil minus 2pt

{\bf Keywords}\quad \thekeywords

\vglue 10pt plus 0.5fil minus 5pt

{\small  Proposed: \theproposer\hfill Received: \receiveddate\nl
Seconded: \theseconders\hfill 
\ifx\reviseddate\relax                         
Accepted: \accepteddate                        
\else
Revised: \reviseddate                          
\fi}
\eject
}       

\let\maketitlepage\maketitlep
\let\maketitle\maketitlepage

\font\phead=cmsl9 scaled 950
\font=cmsl9 scaled 1050
\font\pnum=cmbx10 scaled 913
\font\lnum=cmbx10 
\font\pfoot=cmsl9 scaled 950
\font=cmsl9 scaled 1050
\ifplaintex
\headline{\vbox to 0pt{\vskip -4.5mm\line{\small\phead\ifnum
\count0=\startpage ISSN 1364-0380 (on line)
1465-3060 (printed) \hfill {\pnum\folio}\else\ifodd\count0\def\\{ }%
\ifx\theshorttitle\relax\thetitle\else\theshorttitle\fi\hfill{\pnum\folio}
\else\def\\{ and }{\pnum\folio}\hfill\ifx\theshortauthors\relax\theauthors
\else\theshortauthors\fi\fi\fi}\vss}}
\footline{\vbox to 0pt{\vglue 0mm\line{\small\pfoot\ifnum\count0=\startpage
\copyright\ \gtp\hfill\else
\gt, Volume \thevolumenumber\ (\thevolumeyear)\hfill\fi}\vss
}}
\else
\makeatletter
\def\@oddhead{{\small\lhead\ifnum\count0=\startpage ISSN 1364-0380 (on line)
1465-3060 (printed) \hfill {\lnum\number\count0}\else\ifodd\count0
\def\\{ }\ifx\theshorttitle\relax \thetitle \else\theshorttitle\fi\hfill
{\lnum\number\count0}\else\def\\{ and }{\lnum\number\count0}
\hfill\ifx\theshortauthors\relax 
\theauthors\else\theshortauthors\fi\fi\fi}}\def\@evenhead{\@oddhead}
\def\@oddfoot{\small\lfoot\ifnum\count0=\startpage\copyright\ \gtp\hfill\else
\gt, Volume \thevolumenumber\ (\thevolumeyear)\hfill\fi}
\def\@evenfoot{\@oddfoot}
\makeatother
\fi

\newwrite\gtoutfile
\long\gdef\makeheadfile{  
{\def\\{, }\def\s{ }
\immediate\openout\gtoutfile head.xxx
\immediate\write\gtoutfile{Proxy-for: \ifx\theasciiauthors\relax
\theauthors\else\theasciiauthors\fi\s<\ifx\theasciiemail\relax\theemail\else\theasciiemail\fi>}
\immediate\write\gtoutfile{\noexpand\\}
\immediate\write\gtoutfile{Authors: \ifx\theasciiauthors\relax
\theauthors\else\theasciiauthors\fi}
{\def\\{ }\immediate\write\gtoutfile{Title: \ifx\theasciititle\relax
\thetitle\else\theasciititle\fi}}
\immediate\write\gtoutfile{Subj-class: GT or SG or MG etc}
\immediate\write\gtoutfile{MSC-class: \theprimaryclass\ifx\thesecondaryclass\relax\else, \thesecondaryclass\fi}
\immediate\write\gtoutfile{Journal-ref: Geom. Topol. \thevolumenumber
(\thevolumeyear) \startpage-\finishpage}
\immediate\write\gtoutfile{Comments: Published by Geometry and Topology at}
\immediate\write\gtoutfile{\s\s http://www.maths.warwick.ac.uk/gt/GTVol\thevolumenumber/paper\thepapernumber.abs.html}
\immediate\write\gtoutfile{\noexpand\\}
\immediate\write\gtoutfile{}
\ifx\theasciiabstract\relax
\immediate\write\gtoutfile{\theabstract}\else
\immediate\write\gtoutfile{\theasciiabstract}\fi
\immediate\write\gtoutfile{}
\immediate\write\gtoutfile{\noexpand\\}
\immediate\write\gtoutfile{}
\immediate\closeout\gtoutfile}}  

\def\maketitlepage{\maketitlep\makeheadfile}
\let\maketitle\maketitlepage

\lognumber{307}
\volumenumber{8}\papernumber{1}\volumeyear{2004}
\pagenumbers{1}{34}
\received{4 February 2003}
\published{18 January 2004}
\accepted{13 January 2004}
\proposed{David Gabai}
\seconded{Martin Bridson, Walter Neumann}

\usepackage{epsf,psfrag}

\def\Box{\hbox{\raise 1pt\hbox{$\sqr74$}}}
\def\S{section }

\newtheorem{theorem}{Theorem}[section]

\newtheorem{lemma}[theorem]{Lemma}

\newtheorem{corollary}[theorem]{Corollary}

\def\startproof{\proof}

\def\H{\mbox{\boldmath{$H$}}}%
\def\R{\mbox{\boldmath{$R$}}}%
\def\Z{\mbox{\boldmath{$Z$}}}%

\begin{document}
\title{Modular circle quotients and PL limit sets}
\author{Richard Evan Schwartz }

\address{Department of Mathematics,
University of Maryland\\College Park, MD 20742, USA}

\email{res@math.umd.edu}

\begin{abstract}
We say that a collection $\Gamma$ of geodesics in the
hyperbolic plane $\H^2$ is a {\it modular pattern\/}
if $\Gamma$ is invariant under the
modular group $PSL_2(\Z)$, if
there are only finitely many
$PSL_2(\Z)$--equivalence classes of geodesics
in $\Gamma$, and if
each geodesic in $\Gamma$
is stabilized by an infinite order subgroup
of $PSL_2(\Z)$.  For instance, any finite union
of closed geodesics on the modular orbifold
$\H^2/PSL_2(\Z)$ lifts to a modular pattern.
Let $S^1$ be the ideal boundary of $\H^2$.
Given two points $p,q \in S^1$ we write
$p \sim q$ if $p$ and $q$ are the endpoints
of a geodesic in $\Gamma$. (In particular
$p \sim p$.)  We will see in
\S 3.2 that $\sim$ is an equivalence relation.  We
let $Q_{\Gamma}=S^1/\sim$ be the quotient
space.  We call $Q_{\Gamma}$ a {\it modular circle
quotient\/}.  In this paper we will give a sense
of what modular circle quotients ``look like'' by
realizing them as limit sets of piecewise-linear
group actions
\end{abstract}

\asciiabstract{We say that a collection Gamma of geodesics in the
hyperbolic plane H^2 is a modular pattern if Gamma is invariant under
the modular group PSL_2(Z), if there are only finitely many
PSL_2(Z)-equivalence classes of geodesics in Gamma, and if each
geodesic in Gamma is stabilized by an infinite order subgroup of
PSL_2(Z).  For instance, any finite union of closed geodesics on the
modular orbifold H^2/PSL_2(Z) lifts to a modular pattern.  Let S^1 be
the ideal boundary of H^2.  Given two points p,q in S^1 we write pq if
p and q are the endpoints of a geodesic in Gamma.  (In particular pp.)
We show that is an equivalence relation.  We let Q_Gamma=S^1/ be the
quotient space.  We call Q_Gamma a modular circle quotient.  In this
paper we will give a sense of what modular circle quotients `look
like' by realizing them as limit sets of piecewise-linear group
actions}



\maketitlepage

\section{Introduction}

In this paper we address the question:
What does a tennis racket look like if it is strung
so tightly that the individual strings collapse into
points?  Rather than consider the
expensive disasters produced
by an actual experiment we will consider
related theoretical objects called
{\it modular circle quotients\/}.

We say that a collection $\Gamma$ of geodesics in the
hyperbolic plane $\H^2$ is a {\it modular pattern\/}
if $\Gamma$ is invariant under the
modular group $PSL_2(\Z)$, if
there are only finitely many
$PSL_2(\Z)$--equivalence classes of geodesics 
in $\Gamma$, and if
each geodesic in $\Gamma$
is stabilized by an infinite order subgroup
of $PSL_2(\Z)$.  For instance, any finite union
of closed geodesics on the modular orbifold
$\H^2/PSL_2(\Z)$ lifts to a modular pattern.
Let $S^1$ be the ideal boundary of $\H^2$.
Given two points $p,q \in S^1$ we write
$p \sim q$ if $p$ and $q$ are the endpoints
of a geodesic in $\Gamma$. (In particular
$p \sim p$.)  We will see in
\S 3.2 that $\sim$ is an equivalence relation.  We
let $Q_{\Gamma}=S^1/\sim$ be the quotient
space.  We call $Q_{\Gamma}$ a {\it modular circle
quotient\/}.

In \cite{S1} we encountered a
certain modular circle quotient as
the limit set of a special
representation of $PSL_2(\Z)$ into
$PU(2,1)$, the group of complex projective
automorphisms of the $3$--sphere $S^3$. 
In \cite{S2} we embedded some related
circle quotients into $S^3$.
In this paper we will
treat all the modular circle quotients,
motivated by the constructions in \cite{S2} but
starting from scratch.  Our aim is to give a sense
of what they look like, by realizing them as
limit sets of piecewise linear group actions.

\subsection{Statement of results}

Let $\Gamma$ be a modular pattern of geodesics.
As we explain in \S 3.1, there
is a well-known tiling of $\H^2$ by ideal triangles which
is invariant under the action of
$PSL_2(\Z)$.
We call this
tiling the {\it modular tiling\/}. 
We define $|\Gamma|$ to be one more than
the number of
geodesics in $\Gamma$ which intersect
a given edge of the modular tiling.
We will see in \S 3.2 that this number is finite.
$|\Gamma|$ is independent of the choice of edge,
by symmetry.

Let $S^n$ be the $n$--sphere.  
Our model for $S^n$ is the double of
an $n$--simplex: $S^n=\Delta_+ \cup \Delta_-$,
where $\Delta_+$ and $\Delta_-$ are two
copies of an $n$--simplex, glued along their
boundaries. A {\it simplex\/} of
$S^n$ is a sub-simplex of either $\Delta_+$
or $\Delta_-$.  Say that a {\it punctured simplex\/}
of $S^n$ is a simplex with its vertices deleted.

A homeomorphism $h$ of
$S^n$ is {\it piecewise linear\/} (or PL) if
there is some triangulation of $S^n$ into
finitely many simplices such that $h$ is affine when
restricted to each simplex in the
triangulation.  The set of PL
homeomorphisms of $S^n$ forms a topological
group ${\rm PL\/}(S^n)$ equipped with the compact--open
topology.
Let $H \subset {\rm PL\/}(S^n)$ be
a subgroup.  A compact
subset $\Lambda \subset S^n$ is the
{\it limit set\/} of $H$ if $H$ acts
{\it properly discontinuously\/} on 
$S^n-\Lambda$ and {\it minimally\/} on $\Lambda$.
Thus $H(x)$ is dense in $\Lambda$ for
every $x \in \Lambda$ and
for any compact $K \subset S^n-\Lambda$,
the set $\{h \in H\mid h(K) \cap K \not = \emptyset\}$
is finite.

\begin{theorem}
\label{one}\label{two}
Let $n=|\Gamma|$.
There is an embedding
$i\co  Q_{\Gamma} \to S^n$ and a
monomorphism $\rho\co  PSL_2(\Z) \to {\rm PL\/}(S^n)$ such that
$i(Q_{\Gamma})$ is the limit set of $\rho(PSL_2(\Z))$.
There is a $\rho(PSL_2(\Z))$--invariant partition
of $S^n-i(Q_{\Gamma})$ into punctured simplices,
the vertices of which are
densely contained in $i(Q_{\Gamma})$.
\end{theorem}

One generalization of a modular pattern is a
$PSL_2(\Z)$--invariant map $f\co  \Gamma \to (0,1]$,
where $\Gamma$ is a modular pattern.
Let $\Box(\Gamma)$ be the space of these maps.
Let $CS(S^n)$ be the space of closed
subsets of $S^n$, given the Hausdorff
topology. (Two subsets are close if each is
contained in a small tubular neighborhood of the other.)
Let ${\rm Mon\/}(PSL_2(\Z),{\rm PL\/}(S^n))$
denote the space of monomorphisms from
$PSL_2(\Z)$ into ${\rm PL\/}(S^n)$ given the
algebraic topology. (Two 
monomorphisms are close if they map the generators
to nearby elements of ${\rm PL\/}(S^n)$.)
The following result
organizes all the modular circle quotients based on
subpatterns $\Gamma'$ of $\Gamma$.

\begin{theorem}
\label{three}
Let $n=|\Gamma|$.   There are
continuous maps
$\Lambda\co  \Box(\Gamma) \to CS(S^n)$ and
$\rho\co  \Box(\Gamma) \to {\rm Mon\/}(PSL_2(\Z),{\rm PL\/}(S^n))$
such that the following is true for all $f \in \Box(\Gamma)$.
The set $\Lambda_f$ is the limit set of
$\rho_f(PSL_2(\Z))$ and
$\Lambda_f$ is homeomorphic to $Q_{\Gamma'}$,
where $\Gamma'=f^{-1}(1)$.
\end{theorem}

The method we use to prove Theorem \ref{three}
is flexible and allows us to make a statement
about more general kinds of circle quotients:

\begin{theorem}
\label{six}
Let $\Gamma'$ be the lift to $H^2$ of an
arbitrary finite union of closed geodesics on a
cusped hyperbolic surface $\Sigma$.
Let $Q_{\Gamma'}$ be the circle quotient
based on $\Gamma'$.  For some $n$ there is
an embedding $i\co  Q_{\Gamma'} \to S^n$ and
a monomorphism $\rho\co  \pi_1(\Sigma) \to 
{\rm PL\/}(S^n)$ such that $i(Q_{\Gamma'})$
is the limit set of $\rho(\pi_1(\Sigma))$.
\end{theorem}

If $\Gamma_1$ and $\Gamma_2$ are 
both modular patterns and
$\Gamma_1 \subset \Gamma_2$
then we have an inclusion
$\Box(\Gamma_1) \hookrightarrow
\partial \Box(\Gamma_2)$.  Assuming
this inclusion implicitly, we say that
a sequence $\{f_m\} \in \Box(\Gamma_2)$
{\it degenerates\/} to $f \in \Box(\Gamma_1)$
if $f_m(\gamma) \to 0$ for all
$\gamma \in \Gamma_2-\Gamma_1$ and
$f_m(\gamma) \to f(\gamma)$ if 
$\gamma \in \Gamma_1$. 
Let $n_j=|\Gamma_j|$ for $j=1,2$.
The $n_2$--simplex has faces
which are $n_1$--simplices.  The doubles
of these $n_1$--simplices are
copies of $S^{n_1}$ contained in $S^{n_2}$.
We call these copies the {\it natural embeddings\/}
of $S^{n_1}$ into $S^{n_2}$.

\begin{theorem}
\label{four}
There is a natural embedding
$i\co  S^{n_1} \hookrightarrow S^{n_2}$
with the following property.
Let $\{f_m\} \in \Box(\Gamma_2)$ be a sequence which degenerates 
to $f \in \Box(\Gamma_1)$.
Then the limit sets
$\Lambda_{f_m}$ converge to $i(\Lambda_f)$.
The restriction of $\rho_{f_m}$ to
$\Lambda_{f_m}$ converges 
to the action of
$i \circ \rho_f \circ i^{-1}$ on
$i(\Lambda_f)$.
\end{theorem}

Theorem \ref{four} covers one case not explicitly
mentioned. In \S \ref{standard} we define a certain
{\it standard representation\/} $\rho_0\co  PSL_2(\Z) \to
{\rm\/}PL(S^1)$.  
If $f_m(\gamma) \to 0$ for all
$\gamma \in \Gamma_2$ then
$\Lambda_{f_m}$ converges to a
naturally embedded circle $i(S^1)$ and
the restriction of
$\rho_{f_m}$ to $\Lambda_{f_m}$ converges to
$i \circ \rho_0 \circ i^{-1}$.

Theorem \ref{four} lets us organize all the
modular circle quotients into a coherent whole.
We define
$\Box=\bigcup_{\Gamma} \Box(\Gamma)$.
Let $S^{\infty}$ be the direct limit
of $S^n$, under the system of natural embeddings.
Let $CS(S^{\infty})$ denote the set of finite
dimensional closed subsets of $S^{\infty}$ equipped
with the Hausdorff topology.  Then Theorem \ref{four}
gives a map
$\Lambda\co  \Box \to CS(S^{\infty})$
such that $\Lambda_f$ is homeomorphic to
$Q_{\Gamma'}$ and contained in a naturally embedded
$|\Gamma|$--dimensional sphere.
Here $\Gamma=f^{-1}((0,1])$ and
$\Gamma'=f^{-1}(1)$. The map $\Lambda$ is
continuous when restricted to each finite
dimensional subspace of $\Box$.

The following construction illustrates the
nature of our results.
List all the vertices $v_0,v_1,v_2...$ of $\Box$ with
$v_0$ being the vertex corresponding to
the $0$--map---i.e.\ the empty pattern.  Let
$\{f_t\mid t \in [0,\infty)\} \subset \Box$
be a continuous path such that
$f|_{[0,n]}$ is
contained in the convex hull of the
vertices $v_0,...,v_n$ and $f_n=v_n$.
Here $n=0,1,2...$.
Then $\Lambda_0$ is just the double of a
line segment.  As $t$ increases
$\Lambda_t$ continuously and endlessly
 crinkles up, assuming the topology of
every modular circle quotient as it goes.

\subsection{Comparisons and speculation}

Here are some possible connections to our results:
\begin{enumerate}
\item Our constructions here are similar in spirit
to our constructions in \cite{S3}, where we
related the modular group to Pappus's theorem and
thereby produced discrete representations of the modular
group into the group of automorphisms of the
real projective plane.
\item Our Theorem \ref{one} seems at least
vaguely related to the general results in \cite{BS} about
embedding the boundaries of hyperbolic groups
into $S^n$.  
\item Some of the combinatorial ideas underlying our
constructions are related to the theory
\cite{S} of coding geodesics on the
modular surface using their cutting sequences.
We can work this out explicitly but don't
do it in this paper.
\item $\Box(\Gamma)$ (with its associated maps)
is like a PL version of Teichmuller space.  
The groups attached to the set $\{f\mid f^{-1}(1)=\emptyset\}$
are like PL quasi-Fuchsian 
groups \cite{B, M} in that
their limit sets are topological circles.
The other groups are like cusp groups on the
boundary of quasi-Fuchsian space.
\end{enumerate}

We elaborate on the fourth item.
$\Box(\Gamma)$ is
both richer and poorer than 
Teichmuller space.  It is
richer because it allows for deformations
which cannot exist in hyperbolic geometry.
There are no nontrivial deformations
of the modular group into ${\rm Isom\/}(\H^3)$
whereas $\dim \Box(\Gamma)$
grows unboundedly with the complexity of
$\Gamma$.  
Indeed, one
possible use of our results
is that they provide a topological model
for degenerating families of representations of
punctured surface groups---i.e.\ finite index
subgroups of the modular group---into a Lie group.
Such families generally are extremely
difficult to construct, let alone study
geometrically. 
Our results give a glimpse of how punctured surface
groups might degenerate when
{\it non-simple\/} closed geodesics on
the surface are pinched.

$\Box(\Gamma)$ is poorer than 
Teichmuller space because it only allows for
degenerations which
occur by pinching
closed geodesics.  We don't get things
like geometrically infinite limits.
It almost goes without saying that
$\Box(\Gamma)$ is geometrically
much poorer than Teichmuller space.
It does not enjoy any of the beautiful
rigid structure \cite{G} of Teichmuller space.

We wonder how our 
results transfer to the more
rigid setting of a
Lie group $G$ acting on a homogeneous space $X$.
We think that it ought to be possible sometimes to 
geometrize our constructions and 
produce representations of $PSL_2(\Z)$ into
$G$ which ``realize'' our PL representations.
The result in
\cite{S1} is an example of this.
On the other hand, we think that
there should be strong restrictions
on the types of circle
quotients for each pair $(G,X)$.
A general {\it restriction result\/}
would provide a new tool in
the study of representations of
surface groups into Lie groups, because it
would help control the possible degenerations.

As far as we know, all the modular
circle quotients are non-planar.  At any rate,
many of them are non-planar and hence
cannot be embedded into $S^2$.
Probably all of the modular circle quotients
can be embedded into $S^3$.
However, such embeddings would probably
be very ``distorted'' in general.  We would like
to quantify this distortion, and relate it to
the complexity of the modular pattern.

We also wonder about how our results work out
for circle quotients based on uniform lattices,
but don't have any idea how to proceed.

\subsection{Some ideas in the proof}

Our main idea is to construct an
object we call a {\it modular block\/} (or
{\it block\/} for short.)
A block is a certain subset $\Omega \subset S^n$
equipped with an order $3$ PL automorphism $\sigma$.
A block is based on
a neat partition of
the $n$--simplex into $3^k-1$ smaller
$n$--simplices.  Here $k=(n+1)/2$, with
$n$ always being odd.
The partition is combinatorially isomorphic
to the $k$--fold join of a triangle (which is
an $n$--sphere) minus
one $n$--simplex.   $\Omega$ is obtained
by deleting $2$ simplices from the partition,
so that $\partial \Omega$ 
consists of $3$ non-disjoint
$n$--simplex boundaries, called
{\it terminals\/}.
The remaining $3^k-3$ simplices partition $\Omega$
and are permuted by $\sigma$.

We will construct an infinite network of blocks
glued together along terminals.  The network
is essentially tree-like but its fine structure is
related to the symbolic coding of
geodesics in $\Gamma$.
It turns out that
$Q_{\Gamma}$ is homeomorphic to the
closure of the block vertices. $PSL_2(\Z)$
is represented as 
a subgroup of the automorphism group of
the network. 
Underlying our block network is a
kind of correspondence between
some hyperbolic geometry objects related
to the modular tiling and some
simplicial objects.  We call this a
{\it simplicial correspondence\/}.
The following table summarizes the correspondence.

\noindent\cl{\def\strt{\vrule width 0pt height 12pt depth 7pt}
\begin{tabular}{|c|c|} \hline
{\bf hyperbolic object\/} & {\bf simplicial object\/}\strt
\\
\hline
\hline
the modular tiling $T$ & modular block network\strt \\
\hline
ideal triangle of $T$ & modular block\strt
\\ \hline
geodesic edge of $T$ & terminal\strt \\
\hline
ideal vertex of $T$; geodesic of $\Gamma$
 & vertex of a block.\strt \\
\hline
circle quotient & closure of the block vertices\strt \\ \hline
\end{tabular}}

Here is a more global point of view.
We can define an abstract simplicial complex $C(\Gamma)$
whose vertices are elements of $\Gamma \cup VT$.
Here $VT$ is the set of ideal vertices of the
modular tiling.
We say that a subset $S \subset \Gamma \cup VT$ is
an {\it abstract simplex\/} if it satisfies the following properties:
\begin{enumerate}
\item There is some ideal triangle $\tau$ of $T$
(not necessarily unique)
such that every $s \in S$ is either an ideal vertex
of $\tau$ or a geodesic of $\Gamma$ which
intersects $\tau$.  We say that 
$\tau$ and $S$ are {\it associated\/}.
\item If $\tau$ is associated to $S$ and
$H_{\tau} \subset PSL_2(\Z)$ is the order--$3$
stabilizer subgroup of $\tau$ then $S$ does not
contain an orbit of $H_{\tau}$.  Moreover,
$S$ is not stabilized by an order $2$
element of $PSL_2(\Z)$.
\end{enumerate}
Evidently $PSL_2(\Z)$ acts on $C(\Gamma)$.
It turns out that the maximal abstract simplices of
$C(\Gamma)$ are $n$--dimensional and that
$C(\Gamma)$ minus the vertices is a combinatorial
$n$--manifold.
There are $3^k-3$ maximal
abstract simplices of $C(\Gamma)$ associated to
each $\tau$.   Our construction gives
an embedding of $C(\Gamma)$ into $S^n$
in such a way that these $3^k-3$
abstract simplices map to the simplices partitioning
the block corresponding to $\tau$.
The embedding conjugates the natural action of
$PSL_2(\Z)$ on $C(\Gamma)$ to a subgroup of the
automorphism group of the block network.
The embedding maps the vertex set of
$C(G)$ to a dense subset of the limit set.

So far we have sketched the proof of
Theorem \ref{one}.   For the remaining results,
our idea is to modify the block network by
a certain $2$--step process.
First, we push the blocks apart from each other
by attaching collar-like sets, which we call
{\it separators\/}, onto the block terminals.
Compare Figure 5.1.
This process allows the topology of the limit
set to vary with the stratum of $\Box(\Gamma)$,
as in Theorem \ref{three}. (Theorem
\ref{six} comes as another application.)
Second, we {\it warp\/} the shapes of the
individual blocks, to allow the
representations associated to $\Box(\Gamma_2)$
to degenerate to the representations associated
to $\Box(\Gamma_1)$, as in Theorem \ref{four}.
The element of $\Box(\Gamma)$ determines both the
shapes of the warped blocks and the shapes of
the separators.

\subsection{Overview of the paper}

We have tried to make this paper completely self-contained.  It only
relies on a few basic ideas from linear algebra, hyperbolic geometry,
and real analysis.  We remark to the interested reader that section 2
and 3 makes for a complete, shorter paper in itself, which proves
Theorem \ref{one}.  Here is a plan of the rest of paper:

{\bf Section 2: Modular blocks}, containing:\qua 
 2.1: The Block Lemma;\qua
 2.2: The details; \qua 
 2.3: $3$--dimensional example.

{\bf Section 3: Theorem \ref{one}}, containing:\qua 
 3.1: The modular tiling; \qua
 3.2: Modular pattern basics; \qua
 3.3: Simplicial correspondences; \qua
 3.4: Embedding the quotient; \qua
 3.5: Block networks; \qua
 3.6: Putting it together.

{\bf Section  4: Modified blocks\/}, containing:\qua 
 4.1: Partial prisms; \qua
 4.2: Separators; \qua
 4.3: Warped blocks; \qua
 4.4: Main construction; \qua
 4.5: Degeneration. 

{\bf Section 5: The rest of the results\/}, containing:\qua 
 5.1: Modified correspondences; \qua
 5.2: Modified block networks; \qua
 5.3: Proof of Theorem \ref{three};\qua
 5.4: Proof of Theorem \ref{six}; \qua
 5.5: Proof of Theorem \ref{four}.\qua

{\bf References}

\rk{Acknowledgements}I would like to thank the IHES for their
hospitality during the writing of an early version of this paper.  

The author is supported by NSF Research Grant DMS-0072607.

\section{Modular blocks}

\subsection{The Block Lemma}

Let $k \geq 2$ be an integer and let $n=2k-1$.
Let $\Delta_0$ be an $n$--simplex. 
We say that a {\it modular block\/} is a
set
\begin{equation}
\label{block}
\Omega={\rm closure\/}(\Delta_{0}-\Delta_{1}-\Delta_{2})
\end{equation}
Where $\Delta_{1}, \Delta_{2} \subset \Delta_{0}$
are $n$--simplices with disjoint interiors and
\begin{enumerate}
\item For any indices $i \not = j$ there
are $k$ vertices common to $\Delta_{i}$ and $\Delta_{j}$,
and $\partial\Delta_{i} \cap \partial \Delta_{j}$ is the convex
hull of these common vertices.
\item There is an order $3$ PL automorphism
$\sigma\co  \Omega \to \Omega$ such
that $\sigma$ is affine on $\partial \Delta_j$, with
orbit $\partial \Delta_0 \to \partial\Delta_1 \to 
\partial \Delta_2 \to \partial \Delta_0$.
\end{enumerate}
We call $\partial \Delta_j$ a {\it terminal\/} of
$\Omega$ for $j=0,1,2$.
We call $\partial \Delta_0$ the {\it outer
terminal\/} and $\partial \Delta_1$ and $\partial \Delta_2$
the {\it inner terminals\/}.

Recall from \S 1.1 that $S^n=\Delta_+ \cup \Delta_-$, where $\Delta_+$
and $\Delta_-$ are two copies of a standard $n$--simplex.
Our model for
$\Delta_{\pm}$ is the convex hull of
the standard basis vectors in $\R^{n+1}=\R^{2k}$.
The goal of this chapter is to prove
\begin{lemma}[Block Lemma]
There exists a
modular block whose
outer terminal is $\partial \Delta_+$.
\end{lemma}

\noindent
\proof[Proof -- modulo some details]
Let $e_1,...,e_k$ be the standard basis vectors in
$\R^k$. For any $r \in \R$, let $r_{(k)}=
(r,...,r) \in \R^k$.
For $j=1,...,k$ we define the following
points of $\Delta_0=\Delta_+$:
\begin{equation}
\label{vectors}
A_j=(e_j,0_{(k)}); \hskip 15 pt
B_j=\frac{1}{2n}(2_{(k)}-e_j,2_{(k)}-e_j);
\hskip 15 pt
C_j=(0_k,e_j);
\hskip 15 pt 
\end{equation}
Let $Y=\{Y_j\}_{j=1}^k$ for
each letter $Y \in \{A,B,C\}$.
Let $\langle \cdot \rangle$ denote the convex hull operation.
Note that $\Delta_0=\langle A \cup C \rangle$.
We define
\begin{equation}
\label{hull}
\Delta_1=\langle A \cup B \rangle ; \hskip 20 pt
\Delta_2=\langle B \cup C  \rangle.
\end{equation}
The sets $A \cup B$ and $B \cup C$ are bases for
$\R^{2k}$.  (See Lemma \ref{basis}.)
Hence $\Delta_1$ and $\Delta_2$ are $n$--simplices.
Define
$u=(1_{(k)},-1_{(k)})$.
We have $B_j \cdot u=0$ for all $j$.  Therefore
$B$ is contained in the hyperplane $u^{\perp}$.
We also have $A_j \cdot u=1$ and $C_j \cdot u=-1$ for
all $j$.  Therefore $u^{\perp}$ separates $A$ from $C$.
Hence $\Delta_{1} \cap \Delta_{2}=\langle B \rangle$.
Since $B \in {\rm int\/}(\Delta_0)$ we have
$\partial \Delta_{0} \cap \partial \Delta_{1}=\langle A \rangle$.
Likewise $\partial \Delta_{0} \cap \partial \Delta_{2}=
\langle C \rangle$.  Thus $\Delta_0$, $\Delta_1$, and
$\Delta_2$ satisfy Condition 1.

Let $X=A \cup B \cup C$.  We define
$\sigma\co  X \to X$ by the action
\begin{equation} A_j \to B_j \to C_j \to A_j;
\hskip 20 pt j=1,...,k. \end{equation}
Equation \ref{hull} implies that
$\sigma$ extends to a self-homeomorphism of
$\partial \Omega$, which is affine on each terminal.
Here $\Omega$ is as in Equation \ref{block}.

Say that a $2k$--element $S \subset X$ is {\it good\/}
if it does not contain any orbits of $\sigma$ and
does not equal $A \cup C$.  There are $3^k-1$
good subsets, two of which are $A \cup B$ and $B \cup C$.
Lemma \ref{basis} below shows that every good set
is a basis.  
We define a {\it good simplex\/} to be a
simplex of the form $\langle S \rangle$,
where $S$ is a good subset.  We can extend the
action of $\sigma$ to any individual good simplex
other than $\Delta_1$ and $\Delta_2$ by the rule
$\sigma (\langle S \rangle) = \langle \sigma(S) \rangle$.

We will show below that $\Delta_0$ is triangulated by the
good simplices.
That is, $\Delta_0=\bigcup_S \langle S \rangle$,
and for all  
good simplices $\langle S_1 \rangle$ and
$\langle S_2 \rangle$, we have
\begin{equation}
\label{key}
\langle S_1 \rangle \cap
\langle S_2 \rangle =
\langle S_1 \cap S_2 \rangle.
\end{equation}
$\Omega$ is triangulated by
the good simplices which are not
$\Delta_1$ or $\Delta_2$, and these
are permuted by $\sigma$.  Equation \ref{key}
implies that all the individual actions of
$\sigma$ on good simplices fit together continuously.
Hence $\sigma\co  \Omega \to \Omega$ satisfies Condition 2.
\endproof

\subsection{The details}

For $b \geq 1$ 
we introduce the $b \times b$ matrix
$\Upsilon_b$ whose $(ij)$th entry is $1$ if $i=j$
and otherwise $2$.
This circulent matrix has
the eigenvalue $2b-1$ with multiplicity $1$ and
the eigenvalue $-1$ with multiplicity $b-1$.
Therefore
\begin{equation}
\label{alternate}
(-1)^{b-1} \det(\Upsilon_b)>0.
\end{equation}
Before proving Lemma \ref{basis} let's consider a
representative example which shows how
$\Upsilon_b$ arises in our calculations.
We take $(k,n)=(3,5)$ and show that the set
$S=\{A_1,B_1,B_2,C_2,A_3,C_3\}$ is a basis 
for $\R^6$.  Let $M$ be the matrix whose
rows are elements of $S$.
If some row has a single
$1$ in the $j$th spot, and $0$s in all other spots,
we change $j$th spots of all the other rows to $0$.
We call this {\it simple
row reduction\/}. 
We use a combination of permutations and
simple row reductions to show that
$\det(M) \not = 0$.
Ignoring the factor of
$\frac{1}{2n}$ in the second and third rows:
$$\left[ \matrix{\underline 1&0&0&0&0&0 \cr
                 1&2&2&1&2&2 \cr
                 2&1&2&2&1&2 \cr
                 0&0&0&0&\underline 1&0 \cr
                 0&0&\underline 1&0&0&0 \cr
                 0&0&0&0&0&\underline 1}\right]
\to
\left[   \matrix{1&0&0&0&0&0 \cr
                 0&2&0&1&0&0 \cr
                 0&1&0&2&0&0 \cr
                 0&0&0&0&1&0 \cr
                 0&0&1&0&0&0 \cr
                 0&0&0&0&0&1}\right]
\to
 \left[  \matrix{1&0&0&0&0&0 \cr
                 0&1&0&0&0&0 \cr
                 0&0&1&0&0&0 \cr
                 0&0&0&1&0&0 \cr
                 0&0&0&0&1&2 \cr
                 0&0&0&0&2&1}\right]
$$
This last matrix obviously has nonzero determinant.
Notice also that $\Upsilon_2$ appears in the
bottom right corner, and $2$ is the cardinality
of $S \cap B$.

\begin{lemma}
\label{basis}
Every good set is a basis of $\R^{2k}$.
\end{lemma}

\startproof
Let $S$ be a good set.
Let $b$ be the cardinality of $S \cap B$.
Using permutations and
simple row reduction we see that
\begin{equation}
\label{nontriv}
\det(M)= s \det \left[ \matrix{ I& 0 \cr
               0 & \frac{1}{2n}\Upsilon_b}\right]
\not = 0.
\end{equation}
Here $s \in \{-1,1\}$ depends on the
number of permutations.
\endproof

\begin{lemma}
\label{intersect}
Let $\langle S_1 \rangle$ be a good simplex.
Each codimension--$1$ face $\langle S' \rangle$ of 
$\langle S_1 \rangle$, which is not a face of $\Delta_0$,
is a face of one other good simplex $\langle S_2 \rangle$.
Equation \ref{key} holds for $\langle S_1 \rangle$
and $\langle S_2 \rangle$.
\end{lemma}

\startproof
We have $S'=S_1-Y_j$ for some
$j \in \{1,...,k\}$ and $Y_j \in \{A_j,B_j,C_j\}$.
Without loss of generality assume $j=1$.
By hypotheses $S' \not \subset A \cup C$. Hence
there is exactly one other way to 
complete $S'$ to a good subset:  Namely, 
$S_2=S' \cup Z_1$, where
$Z_1= \{A_1,B_1,C_1\}-S_1$.
Let $M_Y$ and $M_Z$ denote the matrices
whose rows are the elements of
$S_1$ and $S_2$ respectively.
We require that $Y_1$ and $Z_1$ 
appear in the same rows of $M_Y$ and $M_Z$
respectively and that all other rows coincide.
To verify Equation \ref{key} for
$\langle S_1 \rangle$
and $\langle S_2 \rangle$ it suffices to prove that
$\det(M_Y)/\det(M_Z)<0$.  The idea here is that this
causes $Y_1$ and $Z_1$ to lie on opposite
sides of the hyperplane containing
$\langle S' \rangle$.  By symmetry it suffices to consider
the cases $(Y,Z)=(B,C)$ and $(Y,Z)=(A,C)$.  We will
consider these in turn.

\rk{Case 1}
Let $b$ be the cardinality of $S \cap B$.
Since $S'=S_1-B_1 \not \subset A \cup C$ we have
$b \geq 2$.
Using the operations of
Lemma \ref{basis} we get the
formula in Equation \ref{nontriv} for $\det(M_B)$.
When we perform the same operations on
$M_C$ we get the same matrix as
in Equation \ref{nontriv}, except that
all the $2$'s in one of the rows are changed
to $0$'s.   We can then perform one more
simple row reduction, using this row,
to get 
\begin{equation}
\label{compare2}
\det(M_C) = s \det \left[ \matrix{ I& 0 \cr
               0 & \frac{1}{2n} \Upsilon^{(a)}_b}\right]
\end{equation}
for some $a \in \{1,...,b\}$.
Here $ \Upsilon_b^{(a)}$ is created from $\Upsilon_b$ by
changing the $(aj)$th and $(ja)$th entries of $\Upsilon_b$ from
$2$ to $0$, for all $j \not = a$.  Independent
of $a$ we have
\begin{equation}
\label{crucial}
\det( \Upsilon_b^{(a)})=\det(\Upsilon_{b-1}).
\end{equation}
Equations \ref{alternate}-\ref{crucial} give
$\det(M_B)/\det(M_C)<0$.

\rk{Case 2}
Suppose that $A_1$ and $B_1$ are the
first two rows of $M_A$ and that $C_1$ and
$B_1$ are the first two rows of $M_C$.
Let $M$ be the matrix obtained
by replacing the first row of $M_A$ (or $M_C$) by
$A_1+C_1$.  We have $\det(M)=\det(M_A)+\det(M_C)$.
The first row of $M$ is $(1,0,...,0,1,0,...,0)$.
Using this row for row-reduction we can make
all other rows have zeros in the $1$st and $(k+1)$st
positions.   The last $2k-2$ rows of
$M$ are linearly independent by
Lemma \ref{basis}.
Therefore we can perform a series of row
reductions to change the remaining entries of
the second row of $M$ to $0$s.  
Hence $\det(M)=0$ and $\det(M_A)/\det(M_C)=-1$.
\endproof

\rk{Remark} To be sure we checked all the calculations
entailed by the preceding lemma by computer for the cases
$n=3,5,7,9,11,13$.

\begin{corollary}
\label{union}
$\Delta_0$ is the union of the 
good simplices.
\end{corollary}

\startproof
Let $\Delta'$ be the union of the good
simplices.  $\Delta'$ is closed subset
of $\Delta_0$.  If $\Delta' \not = \Delta_0$ then some
codimension one subset of $\partial \Delta'$ separates
the nonempty ${\rm int\/}(\Delta_0-\Delta')$ from
the nonempty ${\rm int\/}(\Delta')$.  Hence there is
a  good simplex $\langle S_1 \rangle$, 
a codimension $1$ face $\langle S' \rangle$ of
$\langle S_1 \rangle$, and
a point
$x \in {\rm int\/}(\Delta_0) \cap 
{\rm int\/}(\langle S' \rangle) \cap \partial
\Delta'.$
Note that $\langle S' \rangle$ is not a face of $\Delta_0$.
By Lemma \ref{intersect} 
there is a  good simplex
$\langle S_2 \rangle$ which also
has $\langle S' \rangle$ as a face, and
$x \in {\rm int\/}(\langle S_1 \rangle \cup \langle S_2 \rangle)
\subset {\rm int\/}(\Delta')$.
This is a contradiction.
\endproof

\begin{corollary}
\label{key2}
Equation \ref{key} is true for all 
good simplices $\langle S_1 \rangle$
and $\langle S_2 \rangle$ provided
that $\langle S_1 \cap S_2 \rangle$
has codimension less than $3$.
\end{corollary}

\startproof
Lemma \ref{intersect} takes care of the 
codimension $1$ case.  Let
$F=\langle S_1 \cap S_2 \rangle$ have codimension $2$.
We will treat the case when $F$ is
not a face of $\Delta_0$, the other case being
very similar.

There are either $3$ or $4$ ways to complete
$S_1 \cap S_2$ to a good set.
From Lemma \ref{intersect} the
corresponding good simplices just
wind around $F$ in a cyclic fashion.
That is, there is a cyclic ordering to
the simplices, such that consecutive
simplices are as in Lemma \ref{intersect}.
The simplices 
are prevented from winding more
than once around $F$ by the
fact that the total dihedral angle
around $F$ is less than $4 \times \pi=4 \pi$.
The topological situation just described
implies that $\langle S_1 \rangle \cap
\langle S_2 \rangle=F$. 
\endproof

\begin{corollary}
Equation \ref{key} holds for every pair of
 good simplices.
\end{corollary}

\startproof
Let $\Delta=\Delta_0$.
Let $\widehat \Delta$
be the abstract simplicial complex obtained by
gluing together all the  good simplices
along the convex hulls of their common vertices.
We have a tautological map $I\co  \widehat \Delta
\to \Delta$ which maps each abstract version
of a good simplex to its realization as
a subset of $\Delta$.   It suffices to
prove that $I$ is a bijection. 

Let $\widehat \Delta_k$ denote the interior of
the complement
of the codimension--$k$ skeleton of $\widehat \Delta$.
Corollary \ref{key2} implies that $I$ is a local
isometry on $\widehat \Delta_3$.   The point here
is that we just need to look at the links of
interior simplices of codimension $1$ and $2$,
and this is what we have done.

$I$ maps the codimension--$3$ skeleton of $\widehat \Delta$
onto the set of codimension--$3$ faces of
the good simplices.
Since $I$ is onto (by
Corollary \ref{union}) the set
$\Delta_3={\rm int\/}(I(\widehat \Delta_3))$ is obtained
from ${\rm int\/}(\Delta)$ by deleting the
codimension--$3$ faces.  Hence
$\Delta_3$ is open, simply connected and dense.
We can find a local isometry $J$, defined on an
open subset of $\Delta_3$, which is the 
inverse of $I$ where defined.
Since $\Delta_3$ is open and
simply connected, $J$ extends by analytic
continuation to a local isometry on $\Delta_3$.
Since $\Delta_3$ is dense, $J$ extends
to all of $\Delta$.

Since $\widehat \Delta_3$ and $\Delta_3$ are
both simply connected it follows from
analytic continuation that the local isometries
$I \circ J$ and $J \circ I$ are the
identity on $\Delta_3$ and
$\widehat \Delta_3$ respectively.  By
continuity, they are the identity
on $\Delta$ and $\widehat \Delta$
respectively. Hence $I$ is a bijection.
\endproof

\subsection{$3$--dimensional example}

We illustrate our construction by working out the
$3$ dimensional case more explicitly.
To get a $3$--dimensional picture we use the projection
$$
V \to ( V \cdot (1,1,-1,-1),
V \cdot(1,-1,1,-1),V \cdot (1,-1,-1,1) )
$$
Using this projection we have
$$
A_1:\left[\matrix{1\cr 1 \cr 1}\right] \hskip 9 pt
A_2:\left[\matrix{1\cr -1 \cr -1}\right]
\hskip 9 pt
B_1:\left[\matrix{0\cr -\frac{1}{3} \cr 0}\right] \hskip 9 pt
B_2:\left[\matrix{0\cr \frac{1}{3} \cr 0}\right]
\hskip 9 pt
C_1:\left[\matrix{-1\cr 1 \cr -1}\right] \hskip 9 pt
C_2:\left[\matrix{-1\cr -1 \cr 1}\right]
$$
Figure 2.1 shows a projection to the $xy$ plane.
The two tetrahedra on the right are supposed to
fit inside the one on the left, as indicated
by the labels.

\begin{figure}[ht!]\small
\psfrag{A1}{$A_1$}
\psfrag{A2}{$A_2$}
\psfrag{B1}{$B_1$}
\psfrag{B2}{$B_2$}
\psfrag{C1}{$C_1$}
\psfrag{C2}{$C_2$}
\cl{\epsfxsize 3in\epsfbox{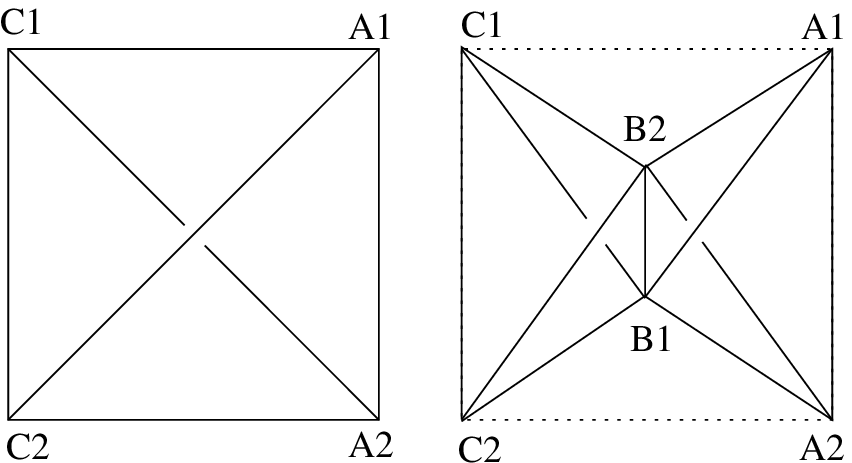}}
\smallskip
\cl{Figure 2.1}
\end{figure}

The $6$ tetrahedra which partition $\Omega$ are
glued together along common faces, in the
following cyclic pattern.
$$
\matrix{&&(A_1B_1C_2B_2)\cr
&\nearrow  && \searrow \cr
(C_1A_1C_2B_2)&&&&(A_1B_1A_2C_2) \cr
&\hskip -15 pt \uparrow && \hskip 15 pt \downarrow \cr
(C_1A_1B_2A_2)&&&&(B_1C_1A_2C_2) \cr
&\nwarrow  && \swarrow  \cr
&&(B_1C_1B_2A_2)}
$$
The action of $\sigma$ translates this cycle of
tetrahedra one third of the way around.
A study of this pattern led us to the
general case.

\section{Theorem \ref{one}}

\subsection{The modular tiling}

We use the disk model of $\H^2$.  By slight
abuse of terminology, we still say that
$PSL_2(\Z)$ acts on this model.
Technically
$PSL_2(\Z)$ acts on the upper half plane model
and a conjugate of $PSL_2(\Z)$ 
acts on the disk model.

$\H^2$ has a canonical (and familiar) tiling $T$ by ideal
triangles which is invariant under the action of
$PSL_2(\Z)$.  We define $T$ by saying that it is the
orbit of an ideal triangle under the group generated
by reflections in its own sides.
See \cite[page 298]{R} for a beautiful picture.
We call $T$ the {\it modular tiling\/}.
Let $VT$, $ET$, and $FT$
respectively denote the set of
ideal vertices, geodesic edges, and ideal triangles
of $T$.  We say that two elements of $ET$ are
{\it touching\/} if they are identical or
share a common endpoint. We say that two
elements of $FT$ are {\it touching\/} if they
are  identical or share a common edge.

$T$ defines an exhaustion of $\H^2$ by ideal polygons.
Let $t_+ \in FT$ be some distinguished ideal triangle.
Let $T_0=\{t_+\}$ and inductively define
$T_{m+1}$ to be those ideal triangles of
$T$ which are touching ideal triangles
of $T_{m}$.
Then $T_m$ is an ideal polygon with
$3 \times 2^m$ sides.  
The {\it combinatorial distance\/}
between $\tau_1,\tau_2 \in FT$ is
the number of edges in $ET$ crossed by
the geodesic segment which connects the centers of
$\tau_1$ and $\tau_2$.  The triangles in
$T_{m+1}-T_m$ are those which are have combinatorial
distance $m$ from $t_+$.

Each $e \in ET$ bounds a unique open halfspace
$h_e$ which is disjoint from the interior of $t_+$.
Given $x \in S^1$ we write
$e|x$ if $x$ is an accumulation point of
$h_e$.  The set of $x$ such that 
$e|x$ is one of the two closed arcs
on $S^1$ determined by the endpoints of $e$.
Each $x \in S^1-VT$ defines a unique
maximal sequence $\{e_m\}$ of edges such that
$e_m|x$ for all $x$ and 
$h_{m+1} \subset h_m$ for all $m$.
We call $\{e_m\}$ the
{\it nesting sequence\/} for $x$.

Using the disk
model of $\H^2$ we can put a metric
on $\H^2 \cup S^1$ which makes it isometric to
a closed Euclidean disk. The next result refers to this metric.

\begin{lemma}
\label{touching}
For any $\eta>0$ there is a $\delta>0$ such that:
If $x_1,x_2 \in S^1$ are
less than $\delta$ apart then there are touching
edges $e_1,e_2 \in ET$ such that
$e_1|x_1$ and $e_2|x_2$.
\end{lemma}

\startproof
There is some $m$ such that all edges of $\partial T_m$,
which is an ideal polygon,
have diameter less than $\eta$.  We take $\delta$ to
be the minimim distance on $S^1$ between vertices of
this ideal polygon.
\endproof

\subsection{Modular pattern basics}

Throughout this chapter $\Gamma$ will denote a
modular pattern of geodesics.

\begin{lemma}
\label{basic}Let $\Gamma$ be a modular pattern of geodesics.
\begin{enumerate}
\item The endpoint of a geodesic of $\Gamma$ belongs
to $S^1-VT$.
\item Two geodesics of $\Gamma$ cannot share an endpoint.
\item Each $e \in ET$ intersects only
finitely many geodesics of $\Gamma$.
\end{enumerate}
\end{lemma}

\startproof
Let $G$ be a finite index torsion-free subgroup of
$PSL_2(\Z)$. 
Then $T/G$ is a tiling of the finite area
surface $\Sigma=\H^2/G$ by ideal
triangles.  Each individual edge of $T$
maps injectively onto an edge of $T/G$.
The quotient
$\Gamma/G$ is a finite union of closed geodesics
on $\Sigma$.
To prove Item 1,
suppose a geodesic $\gamma \in \Gamma$ 
has an endpoint $v \in VT$.
Then $\gamma/G$ 
exits every compact subset of $\Sigma$
as it approaches the cusp point on $\Sigma$
corresponding to $v$.  Closed geodesics
have finite length and hence don't do this.
Item 2 follows from the general fact,
applied to $\gamma_1/G$ and $\gamma_2/G$,
that two closed geodesics on a complete hyperbolic
surface cannot have lifts which share 
exactly one endpoint.
To prove Item 3, note that the set
$(e/G) \cap (\Gamma/G)$ is finite by compactness.
Since the map $e \to e/G$ is injective
the set $e \cap \Gamma$ is also finite.
\endproof

Item 2 above shows that the relation $\sim$ defined
in \S 1 is an equivalence relation:  The transitivity
condition is vacuously satisfied.

\begin{corollary}
\label{basic2}
There is some $m$ such that:
If $\tau_1,\tau_2 \in FT$ have combinatorial
distance at least $m$ then at most one geodesic
of $\Gamma$ intersects both $\tau_1$ and $\tau_2$.
\end{corollary}

\startproof
By symmetry we can take $\tau_1=t_+$, the distinguished
triangle.  If $\{t_m\}$ was a sequence of counterexamples
to this lemma then there would be two geodesics of
$\Gamma$ intersecting both $t_+$ and $t_m$ for all $m$.
Taking a subsequence we can assume that $t_m$ converges
to some $x \in S^1$.
There are only finitely many geodesics which
intersect $t_+$.  Hence, taking another subsequence,
we get the same two geodesics intersecting
$t_m$ for all $m$.  But then $x$ would be an
endpoint to both geodesics, contradicting Item 2 of
Lemma \ref{basic}.
\endproof

\subsection{Simplicial correspondences}

We continue with the notation established above.
For each $e \in ET$ let
$\Gamma_e$ denote the set with the following
description.  An object is an element of $\Gamma_e$
iff it is either an endpoint of $e$ or a
geodesic of $\Gamma$ which crosses $e$.
The cardinality of $\Gamma_e$ is $|\Gamma|+1$,
where $|\Gamma|$ is as in Theorem \ref{one}.
This quantity is finite by Lemma \ref{basic} and
independent of $e$ by
symmetry.  As in the statement of
Theorem \ref{two} we let $n=|\Gamma|$.

Recall from \S 2.1 that
$S^n=\Delta_+ \cup \Delta_-$.
We equip $S^n$ with the piecewise Euclidean metric
inherited from $\Delta_+$ and $\Delta_-$.
As in \S 1.1, a
{\it simplex\/} of $S^n$ is defined
to be a sub-simplex of either
$\Delta_+$ or $\Delta_-$.
If $\Delta$ is an $n$--simplex
of $S^n$
there is a bijective map
from $\Gamma_e$ to $V\Delta$,
the vertex set of $\Delta$,
because the two sets have the same
cardinality.
Compare our table at the end of \S 1.

To each $e \in ET$ we assign a pair
\begin{equation}
\Psi(e)=(\Delta_e,\phi_e),
\end{equation} where $\Delta_e$ is an
$n$--dimensional simplex of $S^n$ and
$\phi_e\co  \Gamma_e \to V\Delta_e$ is a bijection.
When the map $\phi_e$ is not immediately
under discussion we will sometimes abuse
notation and write $\Delta_e=\Psi(e)$.

We say that
$\Psi$ is a
{\it simplicial correspondence\/} for $\Gamma$ if it satisfies
the following $3$ properties:

{\bf Property 1}\qua
For any $\epsilon>0$ there are only finitely
many simplices in the image of $\Psi$ which
have diameter greater than $\epsilon$.

{\bf Property 2}\qua
Given $e_1$ and $e_2$ in $ET$
we let $(\Delta_j,\phi_j)=\Psi(e_j)$.
Suppose $v_j$ is a vertex of $\Delta_j$
for $j=1,2$.  Then $v_1=v_2$ iff
$\phi_1^{-1}(v_1)=\phi_2^{-1}(v_2)$. So,
each vertex in the grand union $\Psi(ET)$
is labelled by a unique element
of $\Gamma \cap VT$.  Compare the
table at the end of \S 1.

{\bf Property 3}\qua
Let $e_1$, $e_2$, $\Delta_1$ and
$\Delta_2$ be as in Property 2.
Let $h_j=h_{e_j}$ for $j=1,2$,
as defined in \S 3.1.
We require that 
${\rm int\/}(\Delta_1) \subset {\rm int\/}(\Delta_2)$
iff $h_1 \subset h_1$ and
${\rm int\/}(\Delta_1) \cap {\rm int\/}(\Delta_2)=\emptyset$
iff $h_1 \cap h_2=\emptyset$.  We also
require that $\partial \Delta_1 \cap \partial \Delta_2$
is the convex hull of their common vertices.
So, the simplices have the same nesting
properties as the open half spaces.

We will construct $\Psi$ in \S 3.6.  First we
want to explore the consequences of its
existence.

\subsection{Embedding the quotient}
\label{embed}

In this section we use $\Psi$ to
define an embedding $i\co  Q_{\Gamma} \to S^n$.
Let $x \in S^1$.
Referring to the notation of \S 3.1, there is a
sequence $\{e_m\}_{m=1}^{\infty} \in ET$ such that
$e_m|x$ for all $m$. (This is true even if
$x \in VT$, but there is not a unique maximal
sequence in this case.)
Let $\Delta_m=\Psi(e_m)$ and
\begin{equation}
\label{nest}
\Psi_{\infty}(x)=\bigcap_{n=1}^{\infty} \Delta_m.
\end{equation}

\begin{lemma}
\label{well}
$\Psi_{\infty}$ is well defined.
\end{lemma}

\startproof
If $x \in VT$ then
$x$ is an endpoint of
$e_m$ for all $m$. Hence $\phi_m(x) \in V\Delta_m$
for all $m$.  Hence
$\Psi_{\infty}(x)=\phi_m(x)$, independent of $m$ 
and the choice of $\{e_m\}$.  
If $x \in S^1-VT$ then any
sequence used to define $\Psi_{\infty}(x)$ is contained in
the nesting sequence for $x$.
The intersection in Equation
\ref{nest} is nested, by Property 3, and is
a single point, by
Property 1. 
\endproof

\begin{lemma}
\label{equiv}
$\Psi_{\infty}$
identifies points on $S^1$ if and only if they are equivalent.
\end{lemma}

\startproof
If $x,x' \in S^1$ are endpoints of
a geodesic $\gamma$ in $\Gamma$ then by
Lemma \ref{basic} we have $x,x' \in S^1-VT$. 
Let $\{e_m\}$ and $\{e_m'\}$ be the nesting
sequences for $x$ and $x'$ respectively.
Then $\gamma$ crosses $e_m$ and $e'_m$ for all $m$ and
$\Delta_m$ and $\Delta'_m$ 
share a vertex for all $m$, by Property 2.  
Hence $\Psi_{\infty}(x)=\Psi_{\infty}(x')$.

If $x,x' \in S^1$ are inequivalent
then there are edges $e,e' \in ET$ such that
\begin{enumerate}
\item $h_e \cap h_{e'}=\emptyset$.
\item $e$ and $e'$ have no vertices in common.
\item $e|x$ and $e|x'$.
\item No geodesic of $\Gamma$ crosses both
$e$ and $e'$.
\end{enumerate}
If this was false then we could take a limit
of a sequence of counterexamples and produce
a geodesic of $\Gamma$ whose endpoints were
$x$ and $x'$.

Now $\Gamma_{e} \cap \Gamma_{e'}=\emptyset$ by
Items 2 and 4.
By Property 2, the simplices
$\Delta_e$ and $\Delta_{e'}$ corresponding
to $e$ and $e'$ have no vertices in common.
Hence $\Delta_e \cap \Delta_{e'}=\emptyset$ by
Property 3 and Item 1.  From the definition of
$\Psi$ and Item 3 we have
$\Psi_{\infty}(x) \in \Delta_e$ and
$\Psi_{\infty}(x') \in \Delta_{e'}$.  Hence
$\Psi_{\infty}(x) \not = \Psi_{\infty}(x')$.
\endproof

\begin{lemma}
\label{cont}
$\Psi_{\infty}$ is continuous.
\end{lemma}

\startproof
Let $\|\cdot \|$ denote the diameter in the
piecewise Euclidean metric on $S^n$ 
and also the Euclidean diameter on
$\H^2 \cup S^1$.
Let $\epsilon>0$ be given.  By Property 1 there is
some $\eta>0$ such that:
If $e \in ET$ satisfies
$\|e\|<\eta$ then
$\|\Delta_{e}\|<\epsilon/2$.
Here $\Delta_e=\Psi(e)$.
Let $\delta$ be as in Lemma \ref{touching}.
If ${\rm dist\/}(x_1,x_2)<\delta$ then
there are touching $e_1,e_2 \in ET$
such that
$e_j|x_j$ and $\|e_j\|<\eta$ for $j=1,2$.
But then
$\Psi_{\infty}(x_1)$ and $\Psi_{\infty}(x_2)$
are contained in simplices
$\Delta_1$ and $\Delta_2$
which by Property 2 share at least one vertex.
Moreover
$\|\Delta_j\|<\epsilon/2$.
Hence ${\rm dist\/}(\Psi_{\infty}(x_1)-\Psi_{\infty}(x_2))<\epsilon$.
\endproof

Define
\begin{equation}
\Lambda=\Psi_{\infty}(S^1).
\end{equation}
Combining the last two results we see that
$\Psi_{\infty}$ factors through a continuous bijection
$i\co  Q_{\Gamma} \to \Lambda$.  
A continuous bijection from a compact space to a Hausdorff space
is a homeomorphism.  Thus $i$ is 
a homeomorphism.  This is our embedding from Theorem \ref{one}.

Here we give a useful characterization of $\Lambda$.

\begin{lemma}  
\label{nest2}
\begin{equation}
\Lambda=\bigcap_{m=0}^{\infty} \Lambda_m
\hskip 30 pt {\rm where\/} \hskip 30 pt
\Lambda_m=\bigcup_{e \in \partial T_m} \Delta_e.
\end{equation}
\end{lemma}

\startproof
We have $\Lambda \subset \Lambda_m$ by Property 3 and
the definition of $\Psi$.
Any  $y \in \bigcap \Lambda_m$ is contained in an
infinite nested sequence $\{\Delta_m\}$ of
simplices.  By Property 3 the
corresponding sequence $\{e_m\}$ is such that
$e_m|x$ for some $x \in S^1$ and for all $m$.
Thus $\Psi_{\infty}(x)=y$.  Hence
$\bigcap \Lambda_m \subset \Lambda$.
\endproof

\rk{Remark}
As above we set
$\Psi(e)=(\Delta_e,\phi_e)$.  By Property 2
all the local maps
$\{\phi_e|\ e \in ET\}$ piece together
to give a global bijection:
$$\bigg[\Gamma \cup VT=\bigcup_{e \in ET} \Gamma_e\bigg]
\hskip 20 pt
\stackrel{\phi}{\longleftrightarrow}
\hskip 20 pt
\bigg[V\Psi=\bigcup_{e \in ET} V\Delta_e\bigg]$$
If $x \in VT$ then
$\Psi_{\infty}(x)=\phi(x)$.  If
$x$ is an endpoint of a geodesic $\gamma$ of
$\Gamma$ then $\Psi_{\infty}(x)=\phi(\gamma)$.
Therefore $\Lambda$ is the
closure of $V\Psi$.

\subsection{Block networks}
\label{bln}

Let $\Omega_+$
be the modular block from \S 2.
We will only use modular blocks
in $S^n$ which have the following definition:
Let $\Delta \subset S^n$ be a
simplex.  Let $A\co  \Delta_+ \to \Delta$ be
an affine isomorphism.   Our new modular
block is $A(\Omega_+)$.  The outer terminal
is $\partial \Delta$. 
Every two modular blocks
in $S^n$ (that we use) are affinely equivalent.
$A$ maps the
canonical triangulation of $\Omega_+$ to
a canonical triangulation of $A(\Omega_+)$.
Given two modular blocks
$\Omega_1, \Omega_2 \subset S^n$ we let
${\rm Map\/}(\Omega_1,\Omega_2)$ be the
set of triangulation-respecting
PL maps from $\Omega_1$ to
$\Omega_2$.

Given $\tau \in FT$ we define
\begin{equation}
\label{gammatau}
\Gamma_{\tau}=\bigcup_{e \in \partial \tau} \Gamma_e.
\end{equation}
For each edge $e \in ET$, the set
$\Gamma_e$ has $2k$ elements.
Each $\Gamma_e$ shares $k$ of its elements with
another $\Gamma_{e'}$.  Hence $\Gamma_{\tau}$
has $3k$ elements.  An $n$--dimensional modular
block $\Omega$ also has $3k$ vertices.
Let $g \in PSL_2(\Z)$ be the element which
cycles the $3$ edges of $\tau$ in counterclockwise
order.  Let $\sigma$ be the $3$--fold PL symmetry
of $\Omega$.   
We say that a {\it $\tau$--labelling\/} of 
$\Omega$ is a bijection $\phi\co  \Gamma_{\tau} \to V\Omega$
which satisfies 
$\phi \circ g = \sigma \circ \phi$.
Here $V\Omega$ is the vertex set of $\Omega$.

\begin{lemma}
\label{symm}
If $\Omega_1$ and $\Omega_2$ are
two $\tau$--labelled modular blocks
then there is a unique element of
${\rm Map\/}(\Omega_1,\Omega_2)$ which
carries the one labelling to the other.
This element is affine if it matches
up the outer terminals.
\end{lemma}

\startproof
Composing with affine maps we reduce to the case
$\Omega_1=\Omega_2=\Omega_+$.
Note first that ${\rm Map\/}(\Omega_+,\Omega_+)$ is
quite large:  Any permutation
$\pi$ of the the $k$--element set
$A=\{A_j\}_{j=1}^k \subset V\Omega_+$ extends
to an element $I_{\pi} \in {\rm Map\/}(\Omega_+,\Omega_+)$
which is an isometry.  The map $I$ commutes with
$\sigma$, the $3$--fold symmetry of $\Omega_+$, and
(hence) permutes the indices of the vectors
in $B$ and $C$ in the same way it permutes the
indices of the vectors in $A$.

Let $\phi$ and $\phi'$ be two $\tau$--labellings
of $\Omega_+$.
Let $e_0,e_1,e_2$ be the three edges of
$\tau$ and let $\Delta_0,\Delta_1,\Delta_2$
be the simplices associated to $\Omega_+$.
Let $\Gamma_j=\Gamma_{e_j}$.
Let $S=\Gamma_0 \cap \Gamma_1$.
Composing $\phi$ with $\sigma^a$ for some $a \in \{0,1,2\}$
we can assume that $\phi(\Gamma_j)=V\Delta_j=
\phi'(\Gamma_j)$ 
for $j=0,1,2$. 
By symmetry $\phi$ and $\phi'$
are determined by their action on
$S$.  Also $\phi(S)=\phi'(S)$.
Hence there is some permutation $\pi$ of
$A=V\Delta_0 \cap V\Delta_1$
such that $\phi|_S=\pi \circ \phi'|_S$.
Referring to the extension of $\pi$ discussed above, we
have $\phi= I_{\pi} \circ \phi'$.
\endproof

We say that a {\it block network\/} is an
assignment $\tau \to \Omega[\tau]$, for
each $\tau \in FT$.  Here
$\Omega[\tau]$ is a $\tau$--labelled modular
block.  We require that
$\Omega[t_+]=\Omega_+$ and
\begin{enumerate}
\item $\Omega[\tau_1]$ and $\Omega[\tau_2]$ have
disjoint interiors for all $\tau_1 \not = \tau_2$.
\item $\Omega[\tau_1]$ and $\Omega[\tau_2]$ share
a common terminal if $\tau_1$ and $\tau_2$
share an edge.
\item $\Omega[\tau_1]$ and $\Omega[\tau_2]$
share a common vertex $v$ if and only if
the $\tau_1$ label of $v$ coincides with the
$\tau_2$ label of $v$.
\end{enumerate}

\begin{lemma}
\label{net1}
There exists a block network for $\Gamma$.
\end{lemma}

\startproof
We choose an enumeration
$t_+=t_0,t_1,t_2,...$ of the ideal triangles of
$FT$ with the following property:
For any $w \geq 1$, 
each $t_w$ shares an
edge $e$ with some $t_v$ for some $v<w$.
We will define $\Omega_w=\Omega[t_w]$ inductively.
We define $\Omega_0=\Omega_+$ as we must.
We choose some $t_0$--labelling for $\Omega_0$.
Note that $S^n-\Omega_0$ consists
of $3$ disjoint open simplices:
${\rm int\/}(\Delta_-)$ and
${\rm int\/}(\Delta_1)$ and
${\rm int\/}(\Delta_2)$.
We call these simplices {\it holes\/}.
Each edge of $t_0$ corresponds to a hole.

Suppose that $\Omega_0,...,\Omega_{w-1}$
have been defined, and each edge of the polygon
$P_w=\partial(t_0 \cup ... \cup t_{w-1})$ is
associated to an open simplex---i.e.\ a hole---of 
$S^n-\bigcup_{j=1}^{w-1} \Omega_j,$
There is some edge $e$ of $P_w$ which 
bounds $t_w$ and some $v<w$ such that
$t_v$ and $t_w$ share $e$ as an edge.
Let $\Delta$ be the simplex which is the
closure of the corresponding hole
in $S^n$.  Note that $V\Delta$
is already labelled by elements of
$\Gamma_e$.  The labelling comes from
the $t_v$--labelling of $\Omega_v$,
which has $\partial \Delta$ a terminal.
First we choose a
$t_w$--labelling of $\Omega_0$ such
that the outer terminal $\partial \Delta_+$ is labelled
by elements of $\Gamma_e$.  Next we choose
the unique affine isomorphism
$A\co  \Delta_+ \to \Delta$ which matches the
$t_w$--labelling of $\Delta_+$ with the
$t_v$--labelling of $\Delta$.  
We define $\Omega_w=A(\Omega_0)$.
We use $A$ to give $\Omega_w$ a
$t_w$--labelling.   The hole $\Delta$ has
been plugged up but the two inner
terminals of $\Omega_w$ bound two new holes.

Our construction only identifies vertices when
they correspond to the same object of
$\Gamma \cup VT$.  No vertices are identified
by accident because of the way the blocks
are nested.   These same nesting properties
show that all the blocks have disjoint interiors.
Thus we have constructed a block network.
\endproof

\rk{Remark} The axioms for block networks
imply that any block network for $\Gamma$ can
be constructed by our inductive process.
Once we determine the $t_0$--labelling
the rest of the construction is forced. Different
$t_0$--labellings produce geometrically identical networks,
but with the labels permuted.

\subsection{Putting it together}

Let $\Omega[*]$ be our block network.  
We define
$\Psi(e)=(\Delta_e,\phi_e)$, where
$\Delta_e$ is the relevant simplex
of $\Omega[\tau]$ and $\phi_e$ is
the restriction of the $\tau$labelling to
$\Gamma_e$.
Here $\tau$ is
one of the two triangles which
has $e$ as an edge.  From the
block network axioms, either choice of $\tau$
gives the same map.

\begin{lemma}
$\Psi$ is
a simplicial correspondence.
\end{lemma}

\startproof
Properties 2 and 3 are immediate from our construction.
It suffices to check property 1.
Given two edges $e,e' \in ET$ we write
$e \to_1 e'$ if $h_e \subset h_{e'}$ and
if $e$ and $e'$ bound a common ideal triangle
of $T$.  We inductively define
$e \to_{(m+1)} e'$ iff $e \to_m e''$ and
$e'' \to _1 e'$.   We let
$\Delta_e=\Psi(e)$ and $\Delta_{e'}=\Psi(e')$.
By Corollary \ref{basic2}, Property 2, and Property 3, there
is some $m$ such that:
If $e \to_m e'$
then $\partial \Delta_e \cap
\partial \Delta_{e'}$ is at most a
single point.  We fix $m$.

Let $\cal S$ denote
the set of pairs of simplices of the form
$(\Delta_e,\Delta_{e'})$, where $e \to_m e'$.
We say that
two pairs $(\Delta_1,\Delta_1')$ and
$(\Delta_2,\Delta_2')$ in $\cal S$ are
{\it equivalent\/} if there is an affine
map which carries one pair to the other.
Modulo $PSL_2(\Z)$ there are only finitely many
pairs $(e,e')$ with $e \to_m e'$.  Thus,
by the affine naturality in our construction,
$\cal S$ contains finitely many
equivalence classes.

Let $\Delta=\Delta_+$, the standard simplex.
For each equivalence class in $\cal S$ we
define a {\it model pair\/}
$(\Delta,\Delta')$, affine equivalent
to any member of the equivalence class.
If Property 1 fails we can find a
nested sequence
$\Delta_1 \supset \Delta_2 \supset \Delta_3...$
such that $\bigcap \Delta_j$ is more than
a single point.  Such a nested sequence
exists by Property 3.  By taking an evenly
spaced subsequence we can assume that
$(\Delta_j,\Delta_{j+1})$ is a member of
$\cal S$ for all $j$.  At least one
model pair $(\Delta,\Delta')$ is represented
infinitely often.

Let $\lambda'_j$ be the longest edge of $\Delta_{j+1}$.
Let $\lambda_j$ be the longer line segment obtained
by intersecting $\Delta_j$ with the line containing
$\lambda'_j$.
Since $\bigcap \Delta_j$ does not shrink
to a point, ${\rm length\/}(\lambda'_j) \not \to 0$.
Since $\Delta_j$ and $\Delta_{j+1}$
converge to each other as $j \to \infty$,
we have ${\rm length\/}(\lambda'_j) \to {\rm length\/}(\lambda_j)$.
For an infinite collection of indices $j$
there is an affine map $T_j$ which
takes the pair
$(\Delta_j,\Delta_{j+1})$ to the model
pair $(\Delta,\Delta')$. 
Let $\delta_j=T_j(\lambda_j)$ and
$\delta'_j=T_j(\lambda'_j)$.
An affine map respects ratios of
distances on lines.  Hence
${\rm length\/}(\delta'_j) \to {\rm length\/}(\delta_j)$
and $\delta'_j$ is an edge of $\Delta'$ for all $j$.
Everything takes place on the same model so 
the set of possible pairs $(\delta'_j,\delta_j)$ is
finite.  Hence
$\delta_j=\delta'_j$ for large $j$.
Hence $\partial \Delta$ and $\partial \Delta'$
have two distinct points in common, contradicting the choice of $m$.
\endproof

The work in \S 3.4 gives us our embedding.  Now we construct
the representation from Theorem \ref{one}.
Let $g \in PSL_2(\Z)$.   Let
$\tau$ be an ideal triangle of $T$.
Let $\tau'=g(\tau)$.   Let
$\phi\co  \Gamma_{\tau} \to \Omega[\tau]$ be the 
$\tau$--labelling of $\Omega[\tau]$.
Let $\phi'\co  \Gamma_{\tau'} \to \Omega[\tau']$ be the 
$\tau'$--labelling of $\Omega[\tau']$.
By Lemma \ref{symm} there is a unique
$\rho(g,\tau) \in
{\rm Map\/}(\Omega[\tau],\Omega[\tau'])$ such
that 
$\rho(g,\tau) \circ \phi=\phi'$.
In other words $\rho(g,\tau)$ maps the vertex of
$\Omega[\tau]$ labelled by the object $x$ to
the vertex of $\Omega[\tau']$ labelled by
the object $x'=g(x)$.  This
{\it intertwining property\/} implies that
$\rho(g,\tau_1)$ and $\rho(g,\tau_2)$ agree on
any common vertices.
Since $\Omega[\tau_1] \cap \Omega[\tau_2]$ is
contained in a single simplex, on which
both our maps are affine, we see that
$\rho(g,\tau_1)=\rho(g,\tau_2)$ on 
$\Omega[\tau_1] \cap \Omega[\tau_2]$.
Letting $\widehat \Omega=\bigcup \Omega[\tau]$
we see that the $\rho(g,\tau)$ maps
piece together to give a continuous
map $\rho(g)\co  \widehat \Omega \to \widehat \Omega$.
The intertwining property gives
\begin{equation}
\label{rep}
\rho(g_1g_2)=\rho(g_1)\rho(g_2).
\end{equation}
We now show that $\rho(g)$ extends to an
element of PL$(S^n)$.
There is some $m$ such that
$t_+ \subset g(T_m)$.  If $e$ is an
edge of $\tau_1 \in FT-T_{m}$ then
$h_e \cap T_m=\emptyset$.  Therefore
$\tau_1 \not \in g(h_e)$.  Therefore $g(h_e)=h_{g(e)}$.
Therefore $\rho(g,\tau_1)$ identifies the
outer terminals of the two blocks and by
Lemma \ref{symm} is affine.
If $\tau_2$ is an ideal triangle touching 
$\tau_1$ and contained in $h_e$ then
$\Omega[\tau_1]$ and $\Omega[\tau_2]$ intersect
along the outer terminal $\partial \Delta$
of $\Omega[\tau_1]$.  
Since two affine maps are determined
by their action on a simplex we see
that $\rho(g,\tau_1)$ and $\rho(g,\tau_2)$
are restrictions of the same affine map.
Repeating this argument with
$\tau_2$ replacing $\tau_1$, etc., we see
inductively that
$\rho(g)$ is affine on all blocks contained
in $\Delta$.  We extend $\rho(g)$ by making
it affine on all of $\Delta$.
Since there are only finitely many
edges of $T_{m+1}$ we see that
the extension of $\rho(g)$ is 
PL on the set $\Lambda_{m+1}$ 
defined in Lemma \ref{nest2}.
There are only finitely many blocks
not contained in $\Lambda_{m+1}$ and
$\rho(g)$ is PL on each one.  In
summary, $\rho(g)$ is a PL map.
Equation \ref{rep} shows that $\rho(g)$ has
the inverse $\rho(g^{-1})$ which is also PL.
Hence $\rho(g) \in {\rm PL\/}(S^n)$.

Equation \ref{rep} says that the map $g \to \rho(g)$
is a homomorphism.  Every $\rho(g)$ acts nontrivially
on some block.  Hence $\rho$ is a monomorphism.
$\Lambda$ is the closure of
the block network vertices.  Hence
$H=\rho(PSL_2(\Z))$ preserves $\Lambda$.    
From the remark at the
end of \S 3.4, the map $\Psi$ conjugates
the minimal action of $PSL_2(\Z)$ on
$\partial \H^2$ to a minimal action of
$H$ on $\Lambda$.
The triangulations of all the
blocks piece together to give a
partition of $S^n-\Lambda$ by punctured
simplices. 
Corollary \ref{basic} and the local finiteness of
the modular tiling imply that our partition
by punctured simplices is locally finite.
$H$ permutes this partition and
hence acts properly discontinuously on
$S^n-\Lambda$.   In short $\Lambda$ is
the limit set of $H$.  Our proof of
Theorem \ref{one} is done.

\section{Modified blocks}

\subsection{Partial prisms}

Let $n=2k-1$ as in previous chapters.
A {\it convex cone\/} in $\R^n$ is a closed convex
subset $C \subset \R^n$, contained in a
halfspace, which is closed under taking
non-negative linear combinations.
$C$ is {\it generated\/} by the set $\Sigma$ if
$C=\{\lambda \Sigma\mid \lambda \geq 0\}$.
We call $C$ a {\it simplex-cone\/}
if $C$ is generated by an $(n-1)$--simplex which 
does not contain $0$.  We also insist that $C$ is
$n$--dimensional.

Let $C$ be a simplex-cone.  Let $H \subset \R^n$
be a codimension $1$ hyperplane which does not
contain $0$.  We say that $H$ {\it cuts\/} $C$ if
$H \cap C$ is an $(n-1)$--simplex.  In this case
$C$ is generated by $H \cap C$.  If $H$ cuts
$C$ we set $\Sigma=C \cap H$ and 
let $[C,H]=\{\lambda \Sigma|\ \lambda \leq 1\}$.
With this definition, $[C,H]$ is an $n$--simplex,
one of whose vertices is $0$, and whose other
vertices are the vertices of $\Sigma$.
We say that a {\it partial prism\/} is a
set isometric to a set of the form
\begin{equation}
\Pi={\rm closure\/}([C,H_1]-[C,H_0]).
\end{equation}
where $H_0$ and $H_1$ cut $C$ and
$[C,H_0] \subset [C,H_1]$.  
We call $C \cap H_0$ the {\it inner boundary\/} of $\Pi$
and we call $C \cap H_1$ the {\it outer boundary\/} of
$\Pi$.  Note---and this is crucial for our constructions---that
the inner and outer boundaries can share vertices in common
or even coincide.  In all cases there is
a canonical bijection between the inner
boundary vertices and the outer boundary vertices:
The matched vertices
lie on the same line through $0$.

$\partial \Pi$ consists of two $(n-1)$--simplices---the inner and
outer boundaries---and some $(n-1)$--dimensional
partial prisms.   This lets us define a
canonical PL involution of $\Pi$
which interchanges
the inner and outer boundaries.  If $n=1$ then
$\Pi$ is an interval or a point, and our involution
reverses the interval or fixes the point depending on the case.  
In general the
PL involution is the cone, to the center of mass of $\Pi$,
of the PL involution which is defined on each
partial prism of $\partial \Pi$ and which
swaps inner and outer boundary components.
We call this map the {\it canonical involution\/}.

We also can define a canonical triangulation of $\Pi$.
If $\Pi$ is a simplex then we use $\Pi$ itself
as the triangulation.  Otherwise we triangulate 
$\partial \Pi$ (by induction) 
and then cone the resulting triangulation to the
center of mass of $\Pi$. We call
this the {\it canonical triangulation\/} of $\Pi$.
The canonical PL involution is affine when restricted
to each simplex in the canonical triangulation.
The important point about the canonical triangulation
is that it has this $2$--fold symmetry.

\subsection{Separators}

Let $n=2k-1$ as above.  Let $\langle \cdot \rangle$ be
the convex hull operation.
We say that a {\it weighted simplex\/} is an
$n$--simplex $\Delta$ together with a map
$S\co  V \Delta \to (0,1]$.  Let
$v_1,...,v_{n+1}$ be the vertices of $\Delta$.
Let $S_i=S(v_i)$.
Let
\begin{equation}
\label{sepp}
v^*_i=S_i v_i + (1-S_i) \beta_S; \hskip 30 pt
\beta_S=\frac{\sum_{i=1}^{n+1} S_i v_i}{\sum_{i=1}^{n+1} S_i}
\in \Delta
\end{equation}
Note that $v^*_i=v_i$ if and only if $S_i=1$.
In all cases $v^*_i$ is contained in the half-open
interval $(\beta_S,v_i]$ which joins $\beta_S$ to $v_i$.
Hence $v^*_1,...,v^*_{n+1}$ are in general
position.  
Finally, we define
\begin{equation}
\label{sepp2}
\Delta_S=\left\langle \bigcup_{i=1}^{n+1} v_i^* \right\rangle;
\hskip 30 pt
[\Delta,S]={\rm closure\/}(\Delta-\Delta_S)
\end{equation}
We call $[\Delta,S]$ a {\it separator\/}.
We call $\partial \Delta$ and $\partial \Delta_S$
respectively the {\it inner\/} and {\it outer
boundaries\/} of $[\Delta,S]$.
For each codimension $1$ face $\Sigma_S$ of $\Delta_S$ there is a
unique codimension $1$ face $\Sigma$ of $\Delta$ such that
$\Sigma_S \subset \langle \beta_S \cup \Sigma \rangle.$
The set
\begin{equation}
\langle \beta_S \cup \Sigma \rangle-
\langle \beta_S \cup \Sigma_S \rangle
\end{equation} is a partial prism.  Therefore
$[\Delta,S]$ is canonically partitioned into $n+1$
partial prisms.

The canonical involutions on the individual prisms
piece together to produce a canonical involution
of $[\Delta,S]$ which swaps
$\partial \Delta$ and
$\partial \Delta_S$.  This involution is affine on each of
$\partial \Delta$ and $\partial \Delta_S$.  The canonical
triangulations of the individual prisms piece
together to give a {\it canonical triangulation\/} of
the separator. The canonical involution of
the separator is affine when restricted to each
simplex in the canonical triangulation.
The combinatorial structure of the canonical
triangulation depends entirely on the set 
$S^{-1}(1) \subset V\Delta$.  Any permutation
of $V\Delta$ which respects this set extends
to a PL automorphism of $[\Delta,S]$ which
is affine on each simplex in the
triangulation and also affine on the boundaries.
This observation underlies Lemma \ref{symm1} below.

If $\psi\co  \Delta \to \Delta'$ is an affine
isomorphism which carries $S$ to $S'$ then
$\psi([\Delta,S])=[\Delta',S']$.  Thus the
separator construction is affinely natural.
Also $[\Delta,S]$ varies continuously with $S$.
Finally, a short computation reveals that
$\beta_S$ is the barycenter of $\Delta_S$.

\subsection{Warped blocks}

Say that a {\it weighted block\/} is a block $\Omega$
equipped with a weighting of its vertices,
$S\co V\Omega \to (0,1]$.
Let $(\Omega,S)$ be a weighted block.
Let $\partial \Delta_0$ be the outer terminal of $\Omega$, so that
$\Omega \subset \Delta_0$.  Let $\partial \Delta_1$ and
$\partial \Delta_2$ be the inner terminals of $\Omega$.
Let $v_1,...,v_{n+1}$ be the vertices of
$\Delta_0$.
Every point of
$\Delta_0$ has the form
\begin{equation}
x=\sum_{i=1}^{n+1} \lambda_i v_i; \hskip 30 pt
{\rm where\/} \hskip 30 pt
\sum_{i=1}^{n+1} \lambda_i=1.
\end{equation}
Let $S_i=S(v_i)$.
We define
\begin{equation}
\label{warp}
P_S(x)=\frac{\sum_{i=1}^{n+1} S_i \lambda_i v_i}
{\sum_{i=1}^{n+1} S_i \lambda_i} \in \Delta_0.
\end{equation}
The map $P_S$ is not a linear map.  However it is
a projective automorphism of $\Delta_0$.
In particular $P_S$ permutes the set of
simplices contained in $\Delta_0$.
We define
\begin{equation}
\Omega_S=P_S(\Omega).
\end{equation}
The terminals of $\Omega_S$ are the
three simplex boundaries
\begin{equation}
\partial \Delta_0; \hskip 20 pt
\Delta_{1,S}=P_S(\partial \Delta_1); \hskip 20 pt
\Delta_{2,S}=P_S(\partial \Delta_2).
\end{equation}
$P_S$ maps the triangulation of $\Omega$ to a
combinatorially equivalent triangulation
of $\Omega_S$.   Thus $\Omega_S$ has a
canonical triangulation.  There is a
canonical PL homeomorphism
$W_S\co  \Omega \to \Omega_S$ which is affine
on each simplex of the triangulations.
$W_S$ conjugates the $3$--fold  PL symmetry $\sigma$ of
$\Omega$ to a $3$--fold PL symmetry $\sigma_S$ of
$\Omega_S$.  By construction $\sigma_S$ is
affine when restricted to each of the terminals
of $\Omega_S$.  

We call $\Omega_S$ a
{\it warped block\/}.
The map $W_S$ sets up a canonical bijection
between $V\Omega$ and $V\Omega_S$.
In this way we transfer the map
$S\co  V\Omega \to (0,1]$ to a map
$S\co  V\Omega_S \to (0,1]$.  In other
words the vertices of $\Omega_S$ are
naturally weighted.  Note that $S$ restricts
to give a weighting to each of the terminals of
$\Omega_S$.  For instance, we have the restriction map
$S\co  V\Delta_{1,S} \to (0,1]$.  

If $(\Omega_1,S_1)$ and $(\Omega_2,S_2)$ are weighted
blocks and $T\co  \Omega_1 \to \Omega_2$ is an affine
map such that $S_1=S_2 \circ T$ then $T(\Omega_{S_1})=
\Omega_{S_2}$.  This follows from the fact that
$T$ conjugates $P_{S_1}$ to $P_{S_2}$, as can be
seen from Equation \ref{warp}.  Our warping
construction is affinely natural even though the map
$P_S$ is not itself affine.

\subsection{Main constructions}

\rk{General modified blocks}
Let $(\Omega,S)$ be as above.
Let $\Delta_S$ be as in Equation \ref{sepp2}.
Let $T_S$ be the affine
map which carries $\Delta_S$ to $\Delta$ in
such a way that $T_S(v^*_i)=v_i$ for all $i$.
Note that $T_S([\Delta,S])$ is a separator whose
inner boundary is $\Delta$.
We define 
\begin{equation}
\label{general}
[\Omega,S]=\Omega_S \cup T_S([\Delta_0,S]) \cup
[\Delta_{1,S},S] \cup [\Delta_{2,S},S].
\end{equation}
We have attached one separator to each terminal
of $\Omega_S$.  We call $[\Omega,S]$ a
{\it general modified block\/}.  We call $\Omega_S$ the
{\it core\/} of $[\Omega,S]$.  Note that $[\Omega,S]$
again has three terminals;  these are the free
boundaries of the attached separators.  The outer
terminal is $T_S(\partial \Delta_0)$.

\rk{Remarks}

(i)\qua Note that 
$\partial \Delta_0$ is the inner boundary of
$T_S([\Delta_0,S])$ whereas
$\partial \Delta_j$ is the outer boundary of
$[\Delta_j,S]$.  From a PL standpoint this
asymmetry in our construction disappears:
Each separator has its canonical involution
which turns it inside out. 

(ii)\qua $[\Omega,S]$ is not necessarily a subset
of $S^n$.  The problem is that the outer boundary of
$T_S([\Delta_0,S])$ might be so large that it is
not contained in one of the two unit simplices
comprising $S^n$.  This difficulty will be
handled in \S 5 in an automatic way.  Our
construction will only use modified blocks
which are contained in $S^n$. 

(iii)\qua Given our definitions in Equations \ref{sepp2}
and \ref{general}, each
vertex $v \in V\Omega_S$
corresponds to two vertices $v_1,v_2 \in V[\Omega,S]$
and we have $v_1=v_2$ if and only if $S(v)=1$.
The same remarks apply to $]\Omega_+,S[$ below.

\rk{Special modified blocks}
Suppose now that $\Omega_+$ is the modular
block constructed in \S 2.  Then
$\Omega_+={\rm closure\/}(S^n-\Delta_--\Delta_1-\Delta_2).$
This follows from the fact that the
outer terminal of $\Omega_+$ is $\partial \Delta_+$, and
$S^n=\Delta_+ \cup \Delta_-$.  The weighting
$S$ gives a map $S\co  V\Delta_- \to (0,1]$ as
well as the maps $S\co  V\Delta_{j,S} \to (0,1]$ for $j=1,2$.
We define
\begin{equation}
\label{specialterminal}
]\Omega_+,S[=(\Omega_+)_S \cup [\Delta_-,S] \cup [\Delta_{1,S},S]
\cup [\Delta_{2,S},S].
\end{equation}
The first separator is contained in
$\Delta_-$.  We call $]\Omega_+,S[$ a {\it special
modified block\/}.  We call $(\Omega_+)_S$ the
{\it core\/} of $]\Omega_+,S[$.  The free
boundaries of the separators are the terminals.

\subsection{Degeneration}

Let $n_1<n_2$ be two integers.
Let $\{\Delta_m\}$ denote a sequence of
$n_2$--simplices.  Let $\Delta'$ be an $n_1$--simplex.
We say that $\Delta_m$
{\it converges barycentrically\/} to $\Delta'$
if some collection of $n_1$ vertices of $\Delta_m$ converges
to the vertices of $\Delta'$ as $m \to \infty$
and the remaining
vertices of $\Delta_m$ converge to the
barycenter of $\Delta'$.  (We shall always have a
consistent labelling of the vertices.)
Referring to \S 4.2:

\begin{lemma}
\label{bary1}
Let $\Delta$ be an $n_2$--simplex and
let $\Delta'$ be an $n_1$--simplex face of $\Delta$.
Let $S_m\co  V\Delta \to (0,1]$ and
$S'\co  V\Delta' \to (0,1]$ be
such that $S_m(v) \to S'(v)$ if $v \in V\Delta'$
and $S_m(v) \to 0$ otherwise.
$\Delta_{S_m}$ converges
barycentrically to $\Delta'_{S'}$.
\end{lemma}

\startproof
Equation \ref{sepp} extends continuously
to the case when some (but not all) of the $S_i$ are
zero.  Extending $S'$ by the $0$--map
we have $S'=\lim S_m$.  When $S'_i=0$
we have $v_i^*=\beta_{S'}$, the
barycenter of $\Delta'_{S'}$.  
\endproof

Let $\Omega_+^j$ denote the $n_j$--dimensional block
from the Block Lemma.  In general we 
use the notation $X^j$ to refer to an object
associated to $\Omega_+^j$ though sometimes we simplify
the notation.
Referring to Equation \ref{vectors} there is a natural embedding
$i\co  S^{n_1} \to S^{n_2}$
defined by 
$i(A_j^1)=A_j^2$ and $i(C_j^2)=C_j^2$ for
$j=1,...,k_1$.
Suppose $S^1\co V\Omega_+^1 \to (0,1]$.
We define $(\Omega',S')=i(\Omega_+^1,S_1)$ and
$\Omega=\Omega_+^2$.
Suppose $S_m\co  V\Omega \to (0,1]$ is a sequence
of maps such that
$S_m(v) \to S'(v)$ if $v \in
V\Omega'$ and $S_m(v) \to 0$ otherwise.
Let $]\Omega,S_m[$ and $[\Omega,S_m]$ be the
special and general modified blocks based on
$(\Omega,S_m)$.
Say that a {\it filled-in terminal\/} of a modified
block is a simplex bounded by a terminal.
The following result is the key to Theorem \ref{four}.

\begin{lemma}
\label{deg1}
The filled-in terminals of 
$]\Omega,S_m[$ converge barycentrically to the
filled-in terminals of $]\Omega',S'[$.  
\end{lemma}

\startproof
By Equations \ref{sepp2} and
\ref{specialterminal} the filled-in
terminals of $]\Omega,S_m[$ are 
$(\Delta_-^2)_{S_m}$ and
$(\Delta^2_{j,S_m})_{S_m}$.
Let $\Delta'=i(\Delta_+^1)$ and
$\Delta'_-=i(\Delta_-^1)$ and
$\Delta'_j=i(\Delta_j^1)$.
The filled-in terminals of $]\Omega',S'[$ are
$(\Delta_-')_{S'}$ and
$(\Delta'_{j,S'})_{S'}$. In all cases, $j \in \{1,2\}$.
Now, $S_m$, $\Delta_+^2$, $S'$ and $\Delta'$ are as
in Lemma \ref{bary1}.  Hence
$(\Delta_-^2)_{S_m}$ converges
barycentrically to $(\Delta_-')_{S'}$.
A direct calculation (which we did numerically on
examples to be sure) shows that the first
$2k_1$ vertices of
$\Delta_{j,S_m}^2=P_{S_m}(\Delta_j^2)$
converge to the vertices of 
$\Delta'_{j,S'}=P_{S'}(\Delta'_j)$.
Lemma \ref{bary1}
finishes the proof in this case.
\endproof

\section{The rest of the results}

\subsection{Modified correspondences}

Suppose that
$\Gamma$ is a modular pattern and $\Gamma' \subset
\Gamma$ is a modular sub-pattern.
We define $\Gamma'_e$ just as we defined
$\Gamma_e$.  We have $\Gamma'_e \subset \Gamma_e$
for all $e \in ET$. 
To each $e \in ET$ we assign a pair
$\Psi'(e)=(\Delta_e,\phi_e)$,
where $\Delta_e$ is an
$n$--dimensional simplex of $S^n$ and
$\phi_e\co  \Gamma_e \to V\Delta_e$ is a bijection.
This is as in \S 3.3.
We say that $\Psi'$ is a
{\it modified simplicial correspondence\/} 
for the pair $(\Gamma',\Gamma)$ if it
satisfies the Properties 1 and 3 for
simplicial correspendences and

\rk{Property 2${}'$}
Given $e_1$ and $e_2$ in $ET$
we let $(\Delta_j,\phi_j)=\Psi'(e_j)$.
Suppose $v_j$ is a vertex of $\Delta_j$
for $j=1,2$.  Then $v_1=v_2$ iff
$\phi_1^{-1}(v_1)=\phi_2^{-1}(v_2)$
{\it and\/} the common object $\phi_j^{-1}(v_j)$
belongs to $\Gamma'_e$.

We define the map $\Psi'_{\infty}\co  S^1 \to S^n$
just as in Equation \ref{nest}.  Lemmas \ref{well},
\ref{cont} and \ref{nest2} work exactly the same way
for $\Psi'_{\infty}$ as they do for $\Psi_{\infty}$.
Property 2${}'$ causes a change in Lemma \ref{equiv}.
The same argument in Lemma \ref{equiv} proves
that $\Psi'$ identifies points on $S^1$ if
and only if they are the common endpoints
of a geodesics in $\Gamma'$.  Thus
$\Psi'_{\infty}$ factors through an
embedding of $Q_{\Gamma'}$ into $S^n$.

\rk{Remark}
The remark at the end of \S 3.4 needs to be
modified in the setting here.
Property 2${}'$ gives a bijection
between $\Gamma' \cup VT$ and a certain
subset $V'\Psi' \subset V\Psi'$ of the
block vertices.   
$\Lambda=\Psi'(S^1)$ is the closure
of $V'\Psi'$.

\subsection{Modified block networks}

Each modified block has a canonical triangulation,
obtained from the triangulations on the
core and on the separators.
Suppose that $\Xi_1$ and
$\Xi_2$ are modified blocks with symmetries
$\sigma_1$ and $\sigma_2$.
Let ${\rm Map\/}(\Xi_1,\Xi_2)$ be
the set of triangulation-respecting PL maps
from $\Xi_1$ to $\Xi_2$.
Suppose that $\Xi_j$ has a weighted core
$(\Omega_j)_{S_j}$ for $j=1,2$.
We say that the bijection $\psi\co  
V((\Omega_1)_{S_1}) \to V((\Omega_2)_{S_2})$
between the core vertex sets
is a {\it perfect matching\/} if 
$\psi^{-1} \circ \sigma_2 \circ \psi=\sigma_1$ and
$S_2 \circ \psi=S_1$.
In other words, $\psi$ is symmetry-respecting
and weight-respecting. 

\begin{lemma}
\label{symm1}
A perfect matching $\psi$ extends to
an element of ${\rm Map\/}(\Xi_1,\Xi_2)$.
When $\Xi_1$ and $\Xi_2$ are general
modified blocks, this extension is
affine if it matches up the outer terminals.
\end{lemma}

\startproof
The combinatorial structure of the separators
of $\Xi$ only depends on $S$.  The combinatorially
identical triangulations on $\Xi_1$ and
$\Xi_2$ define the extension of $\psi$.
When the outer terminals are matched up,
the extension of $\psi$ to the cores is an
affine map $\widehat \psi$.
This follows from the affine naturality of
the warping process.  It follows from
Equation \ref{sepp} that the map $\widehat \psi$
maps the separators of $\Xi_1$ to
the separators of $\Xi_2$.  This follows
from the affine naturality of the separator
construction.
\endproof

The rest of our constructions
depend on some $f \in \Box(\Gamma)$, which
we fix throughout the discussion.
We say that the elements of
$VT$ have weight $1$.  This convention,
together with $f$, assigns
weights to each element of
$\Gamma_{\tau}$, the set in
Equation \ref{gammatau}.
Let $[\Omega,S]$ be a modified block.
Let $\tau \in FT$ be an ideal triangle.
We say that a $\tau$--{\it labelling\/}
of $[\Omega,S]$ is a bijection
$\phi\co  \Gamma_{\tau} \to V\Omega_S$ such that
$\phi \circ g = \sigma \circ \phi$ {\it and\/}
$S \circ \phi=f$.
All maps above have $\Gamma_{\tau}$ as their range.
As in \S 3.5, the element
$g$ is the order $3$ stabilizer of $\tau$
which cycles the edges counterclockwise.
$\sigma$ is the order $3$ PL symmetry of
$[\Omega,S]$.
So, $\phi$ carries the weights
of $f$ to the weighting of $\Omega_S$.
We make the same definitions 
for $]\Omega_+,S[$.

We have labelled $V\Omega_S$, because
$V[\Omega,S]$, the actual vertex set of $[\Omega,S]$,
generally has more vertices than $\Gamma_{\tau}$
has elements.   Here we describe an
{\it induced labelling\/} of $V[\Omega,S]$.
Let $\partial \Delta$ be
one of the terminals of $[\Omega,S]$.
Then one of the terminals $\partial \Delta'$ of
$\Omega_S$ is such that $\partial \Delta$ and
$\partial \Delta'$ form the boundary of a separator
of $[\Omega,S]$.   As with all separators,
there is a canonical bijection $\beta\co 
V\Delta \to V\Delta'$.  One of the
three edges $e$ bounding $\tau$ is such
that $\phi^{-1}(V \Delta)=\Gamma_e$.
We label the vertex $\beta(v) \in V\Delta$
by the pair $(\phi^{-1}(v),e)$. 
In this way, each vertex of $V[\Omega,S]$
is labelled by a pair $(\gamma,e)$, where
$\gamma \in \Gamma_e$ and $e$ is an edge
of $\tau$.  
Given Remark (iii) in \S 4.4, 
and our construction here,
the induced labelling has the property that
$v_1$ is labelled by a pair $(\gamma,e_1)$ and
$v_2$ is labelled by a pair $(\gamma,e_2)$.
Here $\gamma$ is the element of $\Gamma_{\tau}$
which labels $v$. 
We have $v_1=v_2$ if and only if
$f(\gamma)=1$.  We call this the
{\it separation principle\/}.

We define {\it modified block networks\/} 
just as we defined block networks in
\S 3.5, using modified blocks in place
of blocks.  The one twist is that
$\Omega[t_+]$ is a special modified block
and all the other $\Omega[\tau]$ are
general modified blocks. 

\begin{lemma}
\label{net2}
There exists a modified block network for
$f \in \Box(\Gamma)$.
\end{lemma}

\startproof
The proof is essentially the same as the one given
in Lemma \ref{net1}.
Let $t_+=t_0,t_1,t_2,...$ be as in Lemma \ref{net1}.
We need to construct modified blocks
$\Omega_0,\Omega_1,\Omega_2...$,
where $\Omega_j=\Omega[t_j]$.
We set $\Omega_0=]\Omega_+,S[$ as we must.
At the induction step we choose the affine
map $A$ which takes the outer terminal of
$[\Omega_+,S]$ to $\partial \Delta$, the
terminal corresponding to the edge $e$ of $t_v$,
in such a way as to respect the labellings.
\endproof

\rk{Remark} As in \S 3.5 the modified block
network is unique up to the choice of the
$t_0$--labelling.  However, if we base our
construction on some general system of weights
that is not invariant under $PSL_2(\Z)$,
as we do in the proof of Theorem \ref{six} below,
then there are potentially as many different geometric
types of network as their are $G$ equivalence
classes of edges in $ET$.  Here $G$ is the
symmetry group of $f$ (which we will take to
be a finite index subgroup of $PSL_2(\Z)$).

\subsection{Proof of Theorem \ref{three}}

We continue with the notation from above.
We set $\Gamma'=f^{-1}(1)$.  We
use our modified block network to define
a modified correspondence for
$(\Gamma',\Gamma)$:
For each edge $e \in ET$ we define
$\Psi(e)=(\Delta_e,\phi_e)$, where
$\Delta_e$ is the relevant boundary simplex of
$\Omega[\tau]$ and $\phi_e$ is the
labelling of $\Delta_e$ induced by
the $\tau$--labelling of $\Omega[\tau]$.
Here $\tau$ is one of the two ideal
triangles which has $e$ in its
boundary.  Our construction
guarantees that $\Psi'(e)$ is the
same using either choice of $\tau$.
The nesting properties of the modified
block network are the same as for the
original block network.  Hence
$\Psi'$ has property 3.
The same argument as in
\S 3.6 shows that $\Psi'$ has Property 1.

Property 2${}'$ follows from the Separation Principle.
To see how this works, we consider our
construction from a different point of view.
We start with the block
network for $\Gamma$.  We then warp each
block in the network.  This changes the geometry
of the network, but none of its combinatorial
structure.  The warped network still has
property 2. Next, we split apart the 
terminals and insert separators---two per terminal
because the terminal includes into two warped
blocks. (We like to think of this as blowing air
into the terminals.)  Figure 5.1 shows a schematic picture.

The separators have the effect of splitting
apart vertices which are labelled by
geodesics $\gamma$ in $\Gamma$ which have weight less
than $1$.   In the warped block network
$\gamma$ labels a single vertex.  After the
separators are added, there is an infinite
list of vertices associated to $\gamma$.
Each of these vertices has a label of the
form $(\gamma,e)$, where $e$ is an edge of
$ET$ crossed by $\gamma$.  If $\gamma$ has
weight $1$, then all these infinitely
many vertices coalesce into one.
The separators do not affect these weight-$1$ vertices.

\begin{figure}[ht!]\small
\cl{\epsfxsize 3.5in\epsfbox{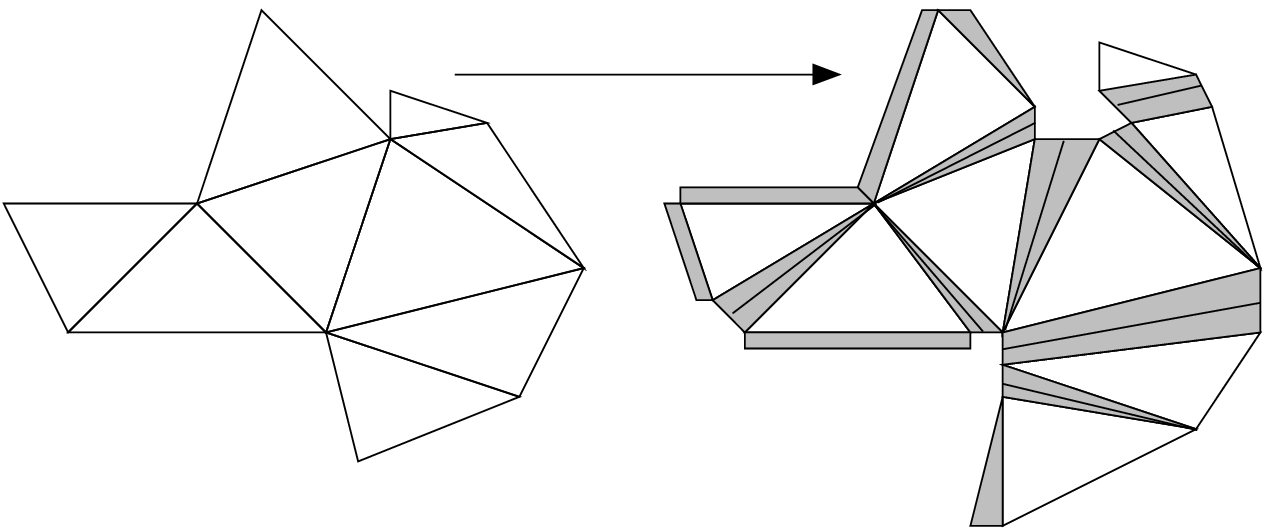}}
\cl{Figure 5.1}
\end{figure}

Thus our modified block network defines a modified
simplicial correspondence $\Psi'$ for $(\Gamma',\Gamma)$.
As in \S 5.1 we have our embedding
$i\co  Q_{\Gamma'} \to S^n$.
We define $\Lambda_f=i(Q_{\Gamma'})$.
The representation $\rho_f$ is constructed
exactly as in \S 3, with Lemma
\ref{symm1} used in place of Lemma \ref{symm}.
The same argument as 
in \S 3 shows that $\Lambda_f$ is the
limit set of $\rho_f$.
The modified blocks and their symmetries are
continuous functions of $f \in \Box(\Gamma)$.
Thus our two maps
$\rho\co  \Box(\Gamma) \to {\rm PL\/}(G,S^n)$ and
$\Lambda\co  \Box(\Gamma) \to [S^n]$ are
continuous maps in the appropriate topologies.

\subsection{Proof of Theorem \ref{six}}

The modular group is hiding behind Theorem
\ref{six}.

\begin{lemma}
\label{modu}
Any cusped finite volume hyperbolic surface
is homeomorphic to a quotient of the form
$\H^2/G$, where $G$ is a finite index
modular subgroup.
\end{lemma}

\startproof
This is a well-known result.  Every cusped
surface has a triangulation into ideal triangles.
Each edge of an ideal triangle has a center point,
the fixed point set of the isometric involution
of the triangle which stabilizes that edge.
We can cut apart our surface and re-glue the
ideal triangles so that the center points
of the edges are matched. This changes
the geometric structure but not the topology.
The resulting
surface then develops into the hyperbolic
plane, onto the modular tiling.  Thus
the new surface, which is homeomorphic
to the original, has the form $\H^2/G$
with $G \subset PSL_2(\Z)$.
\endproof

By Lemma \ref{modu} it suffices to consider the case
of Theorem \ref{six} where $\Sigma=\H^2/G$, so
that $\pi_1(\Sigma)=G$, a finite index modular
subgroup.  
Let $\Gamma$ be the orbit of $\Gamma'$ under
$PSL_2(\Z)$.  Since $G$ has finite index
in $PSL_2(\Z)$, we have that
$\Gamma$ is a modular pattern.
We can define a modified simplicial
correspondence for the pair $(\Gamma',\Gamma)$
even when $\Gamma'$ does not have complete
modular symmetry.  The definitions and
results in \S 5.1 go through word for word.

Let $f\co  \Gamma \to \{\frac{1}{2},1\}$ be the map
defined by the rule $f(\gamma)=1$ if
$\gamma \in \Gamma'$ and
$f(\gamma)=1/2$ if $\gamma \in \Gamma-\Gamma'$.
Even though $f$ is not necessarily
$PSL_2(\Z)$--invariant we can 
define a modified block network for $f$.
This network has $G$--symmetry rather than
$PSL_2(\Z)$--symmetry from the PL standpoint.
The modified block network in turn defines
a modified simplicial
correspondence for $(\Gamma',\Gamma)$.
The same argument as in the proof
of Theorem \ref{three} (which is
just adapted from \S 3) gives a
representation $\rho\co  G \to {\rm PL\/}(S^n)$
which has $i(Q_{\Gamma'})$ as its limit
set.

\subsection{Proof of Theorem \ref{four}}
\label{standard}

Given an affine map $A$ let
$\|A\|=\sup_v \|A(v)\|$ be the operator norm,
with the sup 
being taken over unit vectors.

\begin{lemma} 
\label{lip2}
Let $\Delta'$ be an $n_1$--dimensional face of $\Delta$,
the unit $n_2$--simplex.
Let $\{A_m\}$ be a sequence of affine maps of
$\R^{n_2}$,
with uniformly bounded operator norm, 
such that $A_m|_{\Delta'}$ converges to an
affine injection $A'\co  \Delta' \to \R^{n_2}$.
Let $\{\Xi_m\}$ be a sequence of $n$--simplices
which converge barycentrically to some 
$n'$--simplex $\Xi' \subset \Delta'$.  Then
$A_m(\Xi_m)$ converges barycentrically to
$A'(\Xi')$.
\end{lemma}

\startproof
The first $n_1$ vertices of $\Xi_m$ converge to
the vertices of $\Xi'$.  The bound on the
operator norms guarantees that the first
$n_1$ vertices of $A_m(\Xi_m)$ converge to
$A'(\Xi')$.  The remaining vertices of $\Xi_m$
converge to the barycenter of $\Xi'$.  Again,
the bound on the operator norms guarantees
the images of these remaining vertices under
$A_m$ converge to the barycenter of $A'(\Xi')$.
\endproof

We continue the notation from \S 4.6
and also use the notation from Theorem \ref{four}.
Let $]\Omega,S_m[$ and $[\Omega,S_m]$ be the
special and general modified blocks based on
$\Omega=\ _2\Omega_+$, corresponding to $f_m$.  Thus
$]\Omega,S_m[$ is the zeroth modified
block in the modified block network for
$f_m$ and $[\Omega,S_m]$ is the general
modified block used in the induction step
of Lemma \ref{net2}.   We let $S$ be the
weighting on $_1\Omega_+$ that correponds
to $f \in \Box(\Gamma_1)$. 

For each $m$ we need to choose a $t_+$--labelling
of $]\Omega,S_m[$.  We pick the labelling
so that the vertices in $A' \cup C'$ are
labelled by objects associated to $\Gamma_1$.
We can make the labellings independent of
$m$, since only the weights vary with $m$.
We can choose a $t_+$--labelling of $]_1\Omega,S[$
which is consistent with our $t_+$--labellings
of $]\Omega,S_m[$.  All the same
remarks apply to $[\Omega,S_m]$ and
$[_1\Omega,S]$.  This sets things up so that
$]\Omega,S_m[$ and $[\Omega,S_m]$
are as in Lemma \ref{deg1}.

For each $m$ we have a modified block network
$N_m \subset S^{n_2}$.
We also have a modified block network $N
\subset S^{n_1}$.
Let $N'=i(N)$.  The terminals of
$N_m$ are canonically bijective with the
terminals of $N'$.  Both are indexed by $ET$.

\begin{lemma}
\label{deg2}
Each filled-in terminal of $N_m$ converges barycentrically 
to the corresponding terminal of $N'$ as
$m \to \infty$.  
\end{lemma}

\startproof
Let $t_0,t_1,t_2...$ be as in Lemma \ref{net2}.
For ease of notation we
suppress the dependence on $m$.  Let
$\Omega_j$ be the modified block
associated to $t_j$ when the construction is based on
$f_m \in \Box(\Gamma_2)$.  Let $_1\Omega_j$ be the
modified block
associated to $t_j$ when the construction is based on
$f \in \Box(\Gamma)$.  
Let $\Omega'_j=i(_1\Omega_j)$.  It suffices
to show that the filled-in terminals of $\Omega_j$
converge barycentrically to the filled-in terminals of
$\Omega'_j$.  For $j=0$ this is exactly
Lemma \ref{deg1}.  

Suppose the result is true for $j=1,...,w-1$.  We
consider the case $j=w$.
We adopt the notation from Lemma \ref{net1}
and \ref{net2}.
Thus $D$ is the outer filled-in terminal
of $[\Omega_+,S_m]$ and 
$A\co  D \to \Delta$ is such that
the two filled-in inner terminals of $\Omega_w$ are
$A(\Delta_1)$ and $A(\Delta_2)$.  Here
$\Delta_1$ and $\Delta_2$ are the two inner
filled-in terminals of $[\Omega_+,S_m]$.  By
induction $\Delta$ converges to barycentrically
to one of the inner terminals of $\Omega'_v$.
Thus the outer filled-in terminal of $\Omega_w$ converges
barycentrically to the outer filled-in terminal of
$\Omega'_w$.  We just have to show that
$A(\Delta_1)$ and $A(\Delta_2)$ converge barycentrically
to the inner filled-in terminals of $\Omega'_w$.

Note that
$\Delta_+ \subset D$ and
either $A(D) \subset \Delta_+$ or
$A(D) \subset \Delta_-$.  In either case
$A$ maps the standard unit simplex inside
an isometric copy of itself.  This bounds
$\|A\|$, independent of $m$.
With a view towards using Lemma
\ref{lip2} we let $\Delta'=i(_1\Delta_+)$.
We let $\Xi_m$ be 
$\Delta_1$, the first inner terminal of $\Omega_w$.
We let $\Xi'$ be the first inner terminal $\Delta_1'$
of $\Omega'_w$.   The inner filled-in terminals of
$[\Omega_+,S]$ are the same as two of the terminals
of $]\Omega_+m,S[$.  Therefore, by Lemma \ref{deg1},
we have $\Xi_m \to \Xi'$ barycentrically.
By Lemma \ref{lip2} we see that
$A(\Delta_1) \to A'(\Delta'_1)$ barycentrically.
But $A'(\Delta_1')$ is one of the inner filled-in
terminals of $\Omega_w$.
The same argument works for $\Delta_2$.
This completes the induction step.
\endproof

It follows from Lemma \ref{deg2} 
that the limit sets $\Lambda_{f_m}$
converge to $i(\Lambda)$ and in fact
the maps $S^1 \to \Lambda_{f_m}$
converge pointwise to the map
$S^1 \to \Lambda_f$.  The action of
$\rho_{f_m}$ on $\Lambda_{f_m}$ is determined
by the embedding $S^1 \to \Lambda_{f_m}$ and
by the action of $PSL_2(\Z)$ on
$\partial \H^2$.  Hence
the restriction of $\rho_{f_m}$ to
$\Lambda_{f_m}$ converges to
the restriction of $i \circ \rho_f \circ i^{-1}$
to $\Lambda_f$.  This completes the proof
of Theorem \ref{four}, except in the
case when $\Gamma_1$ is the empty pattern.

\rk{The empty pattern}
To deal with the empty pattern we first
define the {\it standard representation\/}
of $PSL_2(\Z)$ onto $S^1$. 
We think of the unit interval $I_+$ as the
convex hull of the vectors $A=(1,0)$ and
$C=(0,1)$ in $\R^2$.  The midpoint of
$I_+$ is the vector $B=(1/2,1/2)$.
Let $I_{1+}=\langle A \cup B \rangle$
and $I_{2+}=\langle B \cup C \rangle$.
Here $\langle \cdot \rangle$ denotes the
convex hull operation.
Let $I_-$ be another copy of $I+$.
Let $S^1=I_+ \cup I_-$.
Let $\sigma_+ \in {\rm PL\/}(S^1)$ be the
order $3$ element whose action is given by
the orbit $I_- \to I_{1+} \to I_{2+}$ and
$A \to B \to C$.
Let $\iota \in {\rm PL\/}(S^1)$ be the
order $2$ element whose action is given by
the orbit $I_+ \to I_-$ and $A \to C$.
The elements $\sigma_+$ and $\iota$
generate an action on $S^1$ which is
topologically conjugate to the standard
action of $PSL_2(\Z)$ on $\partial \H^2$.
To see this, we identify $A$, $B$, and $C$
with the vertices of the ideal triangle
$t_+$.  Then $\sigma_+$ acts on
$S^1$ just as the order $3$ stabilizer
of $t_+$ acts on $\partial \H^2$ and
$\iota$ acts on $S^1$ just as the
order $2$ stabilizer of one of the edges
of $t_+$ acts on $\partial \H^2$.

When $\Gamma_1$ is the empty pattern we
are dealing with a sequence $\{f_m\}$ which
converges to the $0$--map.  In this case,
the same analysis as that given in
Lemma \ref{deg1} shows that the 
terminals of $]\Omega,S_m[$ converge
to the line segments $I_-$, $I_{1+}$
and $I_{2+}$.  The same argument we give
in Lemma \ref{deg2} then shows that
$\Lambda_{f_m}$ converges to $S^1$ and
the restriction of $\rho_{f_m}$ to
$\Lambda_{f_m}$ converges to the
standard representation.  This completes
our proof of Theorem \ref{four}.

\end{document}